\magnification 1200
\input amstex
\documentstyle{amsppt}
\vsize 7.85in
\hoffset 4.5truemm
\topmatter

\rightheadtext{Convergence of formal CR maps}
\leftheadtext{Jo\"el Merker}
\title On the convergence of S-nondegenerate formal
CR maps between real analytic CR manifolds\endtitle
\author Jo\"el Merker
\endauthor
\address Laboratoire d'Analyse, Topologie et Probabilit\'es,
Centre de Math\'ematiques et d'Informatique, UMR 6632, 39 rue Joliot
Curie, F-13453 Marseille Cedex 13, France. Fax: 00 33 (0)4 91 11 35
52\endaddress
\email merker\@cmi.univ-mrs.fr\endemail
\thanks
\endthanks


\loadeufm

\def\E{{\Bbb E}}
\def\L{{\Bbb L}} 
\def\SS{{\Bbb S}}
\def\N{{\Bbb N}}

\def\R{{\Bbb R}} 
\def\C{{\Bbb C}} 
\def\L{{\Bbb L}} 
\def\V{{\Bbb V}}
\define \dl{[\![}
\define \dr{]\!]}
\def\v{\vert}
\def\dim{\text{\rm dim}}

\def\1{{\text{\bf 1}}} 

\abstract
A new pointwise nondegeneracy condition about generic real analytic
submanifolds in $\C^n$, called {\it S-nondegeneracy}, is
introduced. This condition is intermediate between essential
finiteness and holomorphic nondegeneracy. We prove in this paper that
a formal biholomorphism between S-nondegenerate and minimal real
analytic CR generic manifolds is convergent, generalizing by this
earlier results of Chern and Moser (1974) and of Baouendi, Ebenfelt
and Rothschild (1997-8-9). More generally, we consider S-nondegenerate
formal mappings and establish their convergence. An essential feature
of S-nondegeneracy lies in the unsolvability of the components of the
map in terms of the antiholomorphic components together with their
jets. We conduct convergence results without explicit solvability by
means of repeated applications of Artin's approximation theorem.
\endabstract

\endtopmatter

\document

\head \S 0. Introduction\endhead

In this paper, we study the convergence of formal CR mappings of
smooth complex spaces taking one real analytic submanifold into
another one, following [CM] [BER97] [BER99] and extending the results
therein. We prove here a

\proclaim{Main Theorem}
Any formal invertible CR mapping between minimal Segre nondegenerate
real analytic $(\Cal C^\omega)$ CR manifolds in $\C^n$, $n\geq 2$,
must be convergent.
\endproclaim

\noindent
plus other related generalizations, which cover all analogous results
in the literature.

The study of the convergence of formal holomorphic mappings between
analytic or formal CR objects began with Chern and Moser, who proved
in a celebrated paper ([CM], \S2--3, Theorem 3.5) that the unique
formal transformation taking a Levi-nondegenerate hypersurface in
$\C^n$ ($n\geq 2)$ into normal form is convergent. (Incidentally, the
problem of convergence of formal maps is much deeper and much more
geometric in case one of the object (normal form) is known to be only
formal, and the techniques developed by Chern and Moser cover easily
the case where both are analytic, see also Ebenfelt [E1,2,3,4].)
Later on, in 1997, Baouendi, Ebenfelt and Rothschild proved the main
theorem above assuming that the CR manifolds are finitely
nondegenerate [BER97], based on the implicit function theorem and on
iteration of jet reflection, or essentially finite [BER99], based on
the polynomials identities which had been introduced by Baouendi,
Jacobowitz and Treves in 1985 ([BJT] or [BERbk], Chapter 5). To my
knowledge, these are the only works dealing with the regularity
problem about formal maps between two ${\Cal C}^\omega$ CR manifolds,
although Tr\'epreau told me in 1995 that Treves posed to him in the
eighties the open problem to find a characterization (we provide it in
[Mer99c]).

Quite paradoxically, much attention has been devoted to the
convergence problem for maps between real analytic objects with {\it
CR singularities}, thanks to the influence of Moser and Webster [MW].
In this article, after the works of Bishop [Bi], of Bedford-Gaveau
[BeGa] and of Kenig-Webster's [KW] about the local hull of holomorphy
of a (two-dimensional) surface $S$ in $\C^2$, the authors derived a
complete system of three quantities giving (local) biholomorphic
invariants for $S$ near an (isolated) elliptic complex tangency: such
are biholomorphic to $S_0=\{(z_1,z_2)\in\: y_2=0,\, x_2=z_1\bar z_1 +
(\gamma +\delta x_2^s) (z_1^2+\bar z_1^2)\}$ ($S_0$ is algebraic!),
where $0<\gamma<1/2$ is Bishop's invariant and where $s\in \N$ and
$\delta=\pm 1$, or $s=\infty$ and $\delta=0$. Let us mention further
important results. In 1985, Moser treated the case $\gamma=0$ and
showed that a formal power series change of variables can be found so
that the surface $S$ can be defined by an equation of the form
$\{(z_1,z_2)\: z_2= z_1\bar z_1+ z_1^s+\bar
z_1^s+z_1^{s+1}\varphi(z_1) + \bar z_1^{s+1} \bar \varphi(\bar
z_1)\}$, where $\varphi(z_1)$ is a formal power series in $z_1$, and
where $s$ is a biholomorphic invariant of $S$ at 0. In case
$s=\infty$, Moser observed that $M$ is equivalent to the intersection
$\{y_2=0, x_2=z_1\bar z_1\}$ of the unbounded representation of the
3-sphere with a real hyperplane. In 1995, Huang and Krantz [HK]
completed the study by showing that such elliptic $M$ with $\gamma=0$
is biholomorphically equivalent to $\{(z_1,z_2)\: z_2= z_1\bar z_1+
z_1^s+\bar z_1^s+
\sum_{i+j>s} a_{i,j} z_1^i\bar z_1^j\}$, $\bar
a_{i,j}=a_{j,i}$. Again, as in [CM], one has to deal with
nonnecessarily convergent normal forms. Because of CR singularity,
Moser and Webster could even produce examples of couples of {\it
hyperbolic} $\Cal C^\omega$ surfaces which are formally equivalent but
not biholomorphically equivalent, {\it e.g.} the surface
$\{(z_1,z_2)\: z_2 = z_1\bar z_1+ \gamma \bar z_1^2+\gamma z_1^3 \bar
z_1\}$, $1/2 <\gamma <\infty$, which cannot be biholomorphically
transformed into a real hyperplane, although any surface with $1/2
<\gamma <\infty$ such that the solutions $\mu$ of
$\mu^2+(2-\gamma^{-1})\mu+1=0$ are {\it not} roots of unity can be
{\it formally} transformed into a real hyperplane ([MW] \S5--6). This
divergence is related to the divergence of the normalization of a pair
of involutions $\tau_1$ and $\tau_2$ invariantly attached to the two
two-sheeted projections of the complexification $S^c$ of the surface
$S$ onto the coordinate axes ({\it cf.} also [Ben]), especially to the
{\it elliptic} character of the composition $\varphi=\tau_1\tau_2$
which induces a small divisor problem. Later on, Webster showed [W2]
that each real analytic Lagrangian surface in $\C^2$ with a
nondegenerate complex tangent at $0$ is formally equivalent under
holomorphic {\it symplectic} formal series to the quadric $p=2z\bar
z+\bar z^2$ in $(z,p)\in \C^2$. Again, the non-hyperbolic character
(here, the {\it parabolic} character) of the composition of a similar
pair of involution $\tau_1\tau_2$ enabled Gong [Go2] to show that {\it
generically}, Webster's formal normalization is divergent, following a
suggestion of Moser about divergence of parabolic systems ([Go2],
p.316).

The two special involutions $\tau_1$ and $\tau_2$ attached to the
surface $S$ $(n=2, \hbox{dim}_\R S=2)$ are replaced in the CR
hypersurface (or generic) case $M$ $(n\geq 2, \hbox{dim}_{\R} M
=2n-1$) by the existence of a pair of complexified CR vector fields
$\Cal L$ and $\underline{\Cal L}$ ({\it cf.} [Mer98]) annihilating
formal (anti-)holomorphic mappings, the commutator
$\tau_1\tau_2\tau_1^{-1} \tau_2^{-1}$ (or higher orders commutators)
being the analog of the Levi-form (or of higher order Levi forms, {\it
cf.} Kohn's finite type conditions and Baouendi-Ebenfelt-Rothschild's
finite nondegeneracy) of this pair of vector fields (see in [MW], a
remark p. 262, which almost implicitely suggests the geometric
interpretation of Segre varieties as a double foliation in the
complexified space, as was characterized recently by the author in
[Mer98]). The hypersurface case however simplifies considerably, due
to the constant CR dimension or equivalently, due to the existence of
CR vector fields, which, according to ideas of Sussmann [Su], Treves
[Trv] and Tr\'epreau [Trp], become the mean of {\it propagating
properties of CR functions}. In fact, in this paper and in works
[BER97] [BER99], it is worth noticing that the convergence proofs
reduce to the convergence of formal solutions of ordinary analytic
differential equations with singularities, after remembering that to
any suitably nondegenerate CR hypersurface can be associated a
differentiel equation, {\it cf.} the grounding works of Tresse [Tre]
and Cartan [Ca] (although no published article has yet treated the
correspondence between hypersurfaces and differential equations in a
more degenerate case than the Levi-nondegenerate case; I owe this to
Sukhov). Finally, the natural obstruction to convergence is not due
to a small divisor problem, but to a geometric condition called {\it
holomorphic nondegeneracy} ({\it cf.} [BERbk] [Mer99c]).


Our main intention in this article is to study {\it nonsolvable}
mappings between {\it two} analytic CR manifolds. Thus, the difficulty
does not originate from the possible divergence of some
normalizations, but from the high degeneracies of the {\it Segre
morphism} of the image CR manifold $M'$ ({\it cf.} Introduction of
[Mer99c], and {\it cf.} eq. (1.1.7) below). The finite determination
of formal mappings by their jets at one point, or equivalently what we
call {\it S-solvability} here, has been studied intensively by
Baouendi, Ebenfelt, Rothschild, and Zaitsev and relates strongly to
the solvability (through the usual implicit function theorem only) of
the mapping in terms of the conjugate mapping together with its
conjugate jets. In our analysis, this case appears {\it a posteriori}
to be much simpler ({\it cf.} \S11). We should point out that the main
result in the present paper is new even for invertible mappings, but
leaves open the holomorphically nondegenerate case (except in
codimension 1, {\it cf.} [Mer99c]) which is more general than the
Segre nondegenerate case ({\it cf.} the closing remark in [BER99]).

\head \S 1. Statement of the results\endhead

Let us now explain the words and the concepts in our theorem. Let $h$
be a {\it formal} holomorphic (or CR: this happens to be equivalent)
mapping $(\C^n,p)\to_{\Cal F} (\C^{n'},p')$, {\it i.e.} the components
of $h=(h_j(t-p)=p_j'+\sum_{\gamma\in \N_*^n} h_{j,\gamma} (t-p)^\gamma
)_{1\leq j\leq n'}$, $\N_*^n:=
\N^n\backslash \{0\}$, $h_{j,\gamma}\in \C$, $1\leq j\leq n$,
are {\it formal series} centered in $p$ with respect to the variables
$(t-p)\in \C^n$, with constant term $p'$. We say that $h$ is {\it
invertible} or that $h$ is a {\it formal biholomorphism} or that $h$
{\it has formal rank $n$}, if $n=n'$ and the formal Jacobian of $h$ at
$p$ is invertible. Let $M$ and $M'$ be real analytic CR manifolds in
$\C^n$, $\C^{n'}$, let $p\in M$, $p'\in M'$. We say that the formal
map $h$ {\it maps} $(M,p)$ {\it formally into} $(M',p')$ and write
$h(M,p)\subset_{\Cal F} (M',p')$ or $h\: (M, p) \to_{\Cal F} (M',
p')$, if there exists a $d'\times d$ matrix of formal power series
$\mu(t,\tau)$ such that $\rho'(h(t),\bar{h}(\tau))\equiv \mu(t,\tau)
\rho(t,\tau)$ as formal power series, where $\rho(t,\bar{t})=0$ and
$\rho'(t,\bar{t})=0$ are {\it real analytic} defining equations for
$(M,p)$ and $(M',p')$ respectively. More precisely, our general
assumption throughout this paper will be:

$({\Cal G}{\Cal H})$ The map $h\: (M,p) \to_{\Cal F} (M', p')$ is a
local formal holomorphic map between ${\Cal C}^{\omega}$ CR manifolds
$M \subset \C^n$, $M' \subset \C^{n'}$, $p\in M$ is some point, $p'\in
M'$ is some point, and we assume that $M$ is minimal in the sense of
Tumanov at $p$, or equivalently, of finite type at $p$ in the sense of
Kohn and Bloom-Graham.

Let $m=\hbox{\text{\rm dim}}_{CR} M$, $d=\hbox{\text{\rm codim}}_{\R}
M$, $m'=\hbox{\text{\rm dim}}_{CR} M'$, $d'={\text{\rm codim}}_{\R}
M'$, $m+d=n$, $m'+d'=n'$, $m\geq 1$, $m'\geq 1$, $d\geq 1$, $d'\geq
1$. {\it All our CR manifolds are supposed to be of positive CR
dimension and of positive codimension.}

In suitable coordinates $t$, $t'$, then $p=0$, $f(p)=0$, $M=\{t\in U
\:
\rho(t,\bar{t})=0\}$, $M'=\{t'\in U' \:
\rho'(t',\bar{t}')=0\}$, $U$, $U'$ are small polydiscs centered at the
origin, $\rho_j(t,\bar{t})=\sum_{\mu,\nu\in \N^n} \rho_{j,\mu,
\nu} t^{\mu} \! \bar{t}^{\nu}$, $1\leq j\leq d$,
$\rho_j'(t',\bar{t}')=\sum_{\mu', \nu'\in \N^{n'}} \rho_{j,\mu',
\nu'} {t'}^{\mu'} \! \bar{t'}^{\nu'}$, $1\leq j\leq d'$, are {\it real
analytic}, with $\partial\rho_1 \wedge \cdots \wedge \partial
\rho_d(0)\neq 0$, $\partial\rho_1' \wedge \cdots \wedge \partial
\rho_{d'}'(0)\neq 0$.

By convention, we shall still write sometimes $p$ and $p'$ to denote
the two reference points which are now the origin in the coordinate
systems $t$ and $t'$.

The assumption that $h\: (M,p) \to_{\Cal F} (M', p')$ can be also
interpreted by saying that $\rho'(h(t), \bar{h}(\tau))=0$ when
$\rho(t, \tau)=0$. In particular, $h$ induces a formal holomorphic
map $(S_{\bar{p}}, p) \to_{\Cal F} (S_{\bar{p}'}, p')$, where
$S_{\bar{p}}:=\{t\in U\: \rho(t, 0)=0\}$ is the Segre variety of $M$
at $p=0$ and similarly $S_{\bar{p}'}':=\{ t'\: \rho'(t', 0)=0\}$ is
the Segre variety of $M'$ at $p'=0$. {\it It is this map $h\:
(S_{\bar{p}}, p) \to_{\Cal F} (S_{\bar{p}'}, p')$ which governs the
tangential CR behavior of $f$.}

In $\S1.1$ below, we shall introduce three classes of mappings between
real analytic CR manifolds, one of which is a new class.

\subhead
\S 1.1. Three nondegeneracy conditions on formal CR maps\endsubhead
First, introduce a basis $\underline{\Cal L}_1, \ldots,
\underline{\Cal L}_m$, of $T^{0,1}M$ with complexified coefficients
analytic in $(t, \tau)$.

For instance, after a possible renumbering, we have $\text{\rm det} (
(\frac{\partial \rho_j(0)}{
\partial t_k} )^{1\leq j\leq d}_{m+1 \leq k\leq n})\neq 0$ and
we can thus simply choose the vector fields
$$
\underline{\Cal L}_j= \frac{\partial }{\partial \tau_j} 
-\left(\frac{\partial \rho(t,\tau)}{\partial \tau_j}\right)
\left(\frac{\partial \rho_l (t,\tau)}{
\partial \tau_k}\right)_{1\leq l\leq d;\ m+1 \leq k\leq n}^{-1}
\left(\frac{\partial }{\partial \tau_k}\right)_{m+1\leq k\leq n}.
\tag 1.1.1
$$
Here, $ \left(\frac{\partial }{\partial \tau_k}\right)_{m+1\leq k\leq
n}$ is considered as a $d\times 1$ matrix, $\left(\frac{\partial
\rho(t,\tau)}{\partial \tau_j}\right)$ as a $d\times 1$ matrix and the
$d\times d$ matrix $\left(\frac{\partial \rho_l (t,\tau)}{
\partial \tau_k}\right)_{1\leq l\leq d; \ m+1 \leq k\leq n}$
is invertible, by assumption.

For $\gamma \in \N^m$, denote $|\gamma|:= \gamma_1+\cdots+ |\gamma_m|$
and $\underline{\Cal L}^{\gamma}:=
\underline{\Cal L}_{1}^{\gamma_1}
\cdots \underline{\Cal L}_{m}^{\gamma_m}$. 
Then applying all these derivations to the identity $\rho'(h(t),
\bar{h}(\tau))=0$, it is well-known that one obtains an infinite
family of formal identities
$$
0= \underline{\Cal L}^{\gamma} [\rho'(h(t), \bar{h}(\tau))]:=
R_{\gamma}'(t, \tau, h(t), \nabla^{|\gamma|} \bar{h}(\tau))=0,
\tag 1.1.2
$$
satisfied by $h(t)$ when $(t,\tau)$ satisfy $\rho(t, \tau)=0$. Here,
$\nabla^{|\gamma|} \bar{h}(\tau)$ denotes the
$n'C_{|\gamma|}^{|\gamma|+n}$-tuple of derivatives
$(\partial_{\tau}^{\alpha}\bar{h}(\tau))_{|\alpha|\leq |\gamma|}$ of
$\bar{h}$ with respect to $\tau$ of all orders of lengths $\leq
|\gamma|$, or the $|\gamma|$-jet of $\bar{h}$ at $\tau$, and
$C_{|\gamma|}^{|\gamma|+n}$ denotes Pascal's binomial coefficient
${(\vert \gamma \vert+ n)!\over \vert \gamma \vert ! \ n!}$. Also, one
can easily see that the above term
$R_{\gamma}'=R_{\gamma}'(t,\tau,t',\nabla^{\vert
\gamma\vert})$ denotes here, by its very definition (1.1.2), a
holomorphic mapping from a neighborhood of $0\times 0 \times 0\times
\nabla^{|\gamma|} \bar{h}(0)$ into
$\C^n \times \C^n \times \C^{n'}\times
\C^{n' N_{n,|\gamma|}}$ to $\C^{d'}$, where 
$N_{n, |\gamma|}:= C_{|\gamma|}^{n+|\gamma|}$. Indeed, it is clear
that there exists a holomorphic term $r_\gamma'$ such that we can
write $\underline{\Cal L}^{\gamma}\bar{h}(\tau)= r_{\gamma}' (t, \tau,
\nabla^{|\gamma|}\bar{h}(\tau))$, because coefficients of
$\underline{\Cal L}^{\gamma}$ are analytic in $(t, \tau)$. Denote
$R_{\gamma}'=({R'}_{\gamma}^{l'})_{1\leq l'\leq d'}$ and $R_{\gamma}'=
{R'}_{\gamma}^{l'}(t, \tau, t', \nabla^{|\gamma|} \bar{h}(\tau))=
\underline{\Cal L}^{\gamma} [\rho_{l'}'(h(t), \bar{h}(\tau))]$,
$1\leq l'\leq d'$.

\definition{Definition 1.1.3}
We will say that the formal mapping $h$ is

$\bullet$ {\bf S-solvable} at $p$ if the holomorphic mapping
$$
\C^n \ni t' \mapsto (R_{\gamma}'(0, 0, t', \nabla^{|\gamma|}
\bar{h}(0)))_{|\gamma|\leq \kappa_0} \in 
\C^{d'N_{n,\kappa_0}}
\tag 1.1.4
$$
is an immersion at $0$, for $\kappa_0$ large enough. Then the first
integer $\kappa_0$ for which the mapping in (1.1.4) is an immersion is
in fact a biholomorphic invariant of $h$ under simultaneous changes of
coordinates near $(M, p)$ and near $(M', p')$. Referring to this
integer, we shall shortly say that {\it $h$ is $\kappa_0$-solvable},
in order to means that {\it $h$ is $\kappa_0$-solvable in terms of
$\bar h$ and its jets $\nabla^{\kappa_0} \bar h$} (for the complete
explanation, see in advance Lemma 3.12 in this article). For instance,
it is well-known that $h$ is $l_0'$-solvable in the following
circumstance: when $M'$ is {\it $l_0'$-finitely nondegenerate at $p'$}
in the sense of Baouendi, Ebenfelt and Rothschild [BER97] and $h$ is a
formal {\it CR submersive} map, {\it i.e.} which induces a formal
submersion $(S_{\bar{p}},p)\to_{\Cal F} (S_{\bar{p}'}, p')$ at $p$.

$\bullet$ {\bf S-finite} if the complex analytic variety $\V_p'$
defined by
$$
\V_p':=\{t'\in \C^{n'} \: R_{\gamma}' (0, 0, t', 
\nabla^{|\gamma|} \bar{h}(0))=0,
\ \forall \ \gamma \in \N^m\}
\tag 1.1.5
$$
is zero-dimensional at the point $p'$: $\text{\rm dim}_{\C, p'} \V_p'=
0$ ({\it i.e.} $\text{\rm dim}_{\C, 0} \V_0'= 0$, since $p=0$ and
$p'= 0$ in our coordinates). (The study of $S$-finite CR maps is very
classical and standard, since the work of Baouendi, Jacobowitz and
Treves, see [BJT] [DF88] [BR88] [BR90] [BR95] [BHR96] [CPS98] [BERbk]
[CPS99].)

$\bullet$ {\bf S-nondegenerate} if there exist multiindices $\gamma_1,
\ldots,\gamma_{n'} \in \N^m$ and integers $l_1',\ldots, l_{n'}'$,
$1\leq l_i'\leq d'$, such that
$$
{\text{\rm det}} \left(
\frac{\partial {R'}_{\gamma_j}^{l_j'}}{\partial t_k'}(t, 0, h(t), 
\nabla^{|\gamma_j|} \bar{h}(0)) \right)_{1\leq j,k\leq n'} \not\equiv 0
\tag 1.1.6
$$ \vskip -0.2cm\noindent when $\rho(t, 0)=0$ and where the above
formal series should be interpreted as a formal series expressed in
terms of a local holomorphic coordinate on the Segre variety
$\{\rho(t, 0)=0\}$ passing through $0$. More precisely, as one can
find coordinates $t=(w, z)\in
\C^m\times \C^d$, such that $M$ is given by a $d$-dimensional
vectorial equation in the form $z= \bar{z}+i\bar{\Theta}(\bar{w},
w,\bar{z})$ (see \S2.1, eq. (2.1.1)), the nondegeneracy condition
should be understood as meaning the following:
$$
{\text{\rm det}} \left(
\frac{\partial {R'}_{\gamma_j}^{l_j'}}{\partial t_k'}(w, 
i\bar{\Theta}(0, w, 0), 0, h(w, i\bar{\Theta}(0, w,0)),
\nabla^{|\gamma_j|} \bar{h}(0)) \right)_{1\leq j,k\leq n'} 
\not\equiv_w 0.
\tag 1.1.7
$$
\enddefinition

\vskip -0.2cm\noindent
{\it Remarks.} 1. Using the biholomorphic invariance of the Segre
varieties attached to $M$ and to $M'$, it can be easily shown that
$S$-solvability, $S$-finiteness and $S$-nondegeneracy of a formal CR
map do not depend on the choice of some defining functions
$(\rho_j)_{1\leq j\leq d}$ for $M$ and $(\rho_j')_{1\leq j\leq d'}$
for $M'$, and that these conditions are invariant under simultaneous
biholomorphic changes of coordinates near $M$ and $M'$ which fix $p$
and $p'$.

2. An S-solvable map is clearly S-finite, but there is no general link
between S-finite and S-nondegenerate maps, as shown by simple examples
in \S5 here.

\subhead \S 1.2. The main result\endsubhead
We obtain a general convergence result about formal CR maps between
real analytic CR manifolds satisfying each one of the above three
nondegeneracy conditions. The third is new and constitutes the core of
this article.

\proclaim{Theorem 1.2.1}
Let $h: (M, p)\to_{\Cal F} (M', p')$ be a formal holomorphic map
between real analytic CR generic manifolds and assume that $M$ is
minimal at $p$. If
\roster
\item"(i)" $h$ is S-solvable, or if
\item"(ii)" $h$ is S-finite, or if
\item"(iii)" $h$ is S-nondegenerate,
\endroster
then the power series of the formal mapping $h$ is convergent.
\endproclaim

\remark{Remarks} 
1. An elementary examination of our proof shows that this result
extends immediately to $M'$ being any real analytic set through $p'$,
which is not necessarily smooth nor CR, provided each one of the
nondegeneracy conditions (1.1.4), (1.1.5), or (1.1.6) holds ({\it cf.}
also [CPS99]). However, it seems to be essential in our proof that $M$
is CR generic and minimal and it would be an interesting problem to
search for generalizations of the notion of finite type in the
category of singular real analytic varieties ({\it cf.} [BG]).

2. Although not stated in this form, parts (i) and (ii) of Theorem
1.2.1 were essentially proved in [BER97] and [BER99]
respectively. Furthermore, the versions of (ii) that are proved in
[BER99] followed in fact from [BER97] and earlier techniques developed
by Baouendi and Rothschild in the $\Cal C^\infty$-$\Cal C^\omega$
regularity problem [BJT] [BR88] [BR90].
\endremark

\subhead \S 1.3. Discussion of the proof\endsubhead
Our proof of Theorem 1.2.1 (iii) incorporates two essential
ingredients. As a first ingredient, we shall derive from the {\it
approximation theorem} of Artin ([A]) a beautiful convergence theorem
(Theorem 1.3.2 below). This convergence argument will be applied at
the level of Segre varieties and of subsequent Segre chains. The
approximation theorem states that formal solutions to analytic
equations can be approximated to any order by convergent solutions:

\proclaim{Theorem 1.3.1}
\text{\rm (Artin, [A])}
Let $R(w,y)=0$, $R= (R_1,\ldots,R_J)$, where $w\in \C^n$, $y\in \C^m$,
$R_j\in {\Cal O}_{n+m}=\C\{w,y\}$, be a converging system of \text{\rm
holomorphic} equations. Suppose
$\hat{g}(w)=(\hat{g}_1(w),\ldots,\hat{g}_m(w))$, $\hat{g}_k(w)\in
\C\dl w\dr$, are \text{\rm formal} power series without constant term
which solve $R(w,\hat{g}(w))\equiv_w 0$ in $\C\dl w\dr$. Then for
every integer $N\in \N$, there exists a \text{\rm convergent series
solution} $g(w)= (g_1(w),\ldots,g_m(w))$, {\it i.e.} satisfying
$R(w,g(w))\equiv_w 0$, such that $g(w)\equiv_w \hat{g}(w) \
({\text{\rm mod}} \ \frak{m}(w)^N)$.
\endproclaim

Here, $\frak{m}(w)$ denotes the maximal ideal of the local ring $\C\dl
w\dr$ of formal power series in $w$ and the congruence relation
$g(w)\equiv_w \hat{g}(w) ({\text{\rm mod}} \ \frak{m}(w)^N)$ means
that the coefficients of monomials of total degree $<N$ agree in
$g(w)$ and $\hat{g}(w)$. We denote by $\C\{w\}$ the local ring of
convergent power series in $w$.

\proclaim{Theorem 1.3.2}
Let $R(w,y)=0$, $R=(R_1,\ldots,R_J)$, where $w\in \C^n$, $y\in \C^m$,
$R_j\in {\Cal O}_{n+m}=\C\{w,y\}$ be a system of {\rm holomorphic}
equations. Suppose that
$\hat{g}(w)=(\hat{g}_1(w),\ldots,\hat{g}_m(w))$ $\in\C\dl w\dr^m$ are
{\rm formal} power series without constant term solving
$R(w,\hat{g}(w))\equiv_w 0$ in $\C\dl w\dr$. If $J\geq m$ and if there
exist $j_1,\ldots,j_m$, $1\leq j_1 < j_2 < \cdots < j_m \leq J$ such
that
$$
{\text{\rm det}} \left(
\frac{\partial R_{j_k}}{\partial y_l}(w,\hat{g}(w))
\right)_{1\leq k,l \leq m} \not\equiv_w 0 \ \ {\text{\rm in}} \ \ \C\dl w\dr,
\tag 1.3.3
$$
then $\hat{g}(w)\equiv g(w)\in \C\{w\}$ is convergent.
\endproclaim

This corollary will be of paramount importance in proving Theorem
1.2.1 (iii).

The second main ingredient in our proof of Theorem 1.2.1 (iii) will be
an argument about propagation of analyticity which is closely related
to the recent works of Baouendi, Ebenfelt and Rothschild and which
will be applied here using the formalism that the author have
introduced in [Mer98], a formalism which stems from the local theory
of foliations by flows of vector fields and appears to be canonical
for the following reason ({\it cf.} [Mer98]).

In the extrinsic complexification ${\Cal M}=M^c$ of $M$, the
complexifications of Segre varieties give birth to concatenations of
Segre varieties, called Segre chains in [Mer98] and which do not
coincide {\it exactly} with the so-called Segre sets introduced by
Baouendi, Ebenfelt and Rothschild ({\it cf.} [BER96] [BER97] [Z97]
[BERbk] [BER99]), but coincide up to a change of parametrization. Our
formalism interprets these concatenations of Segre varieties as
partial orbits of the complexified CR vector fields tangent to
$M$. Applying then an iteration processus giving the analyticity of
transversal jets (transversal to subsequent Segre chains), which
follows a each step by subsequent applications of Theorem 1.3.2, we
shall obtain analyticity of $h^c=(h,\bar{h})$ along all Segre
chains. We then conclude by noticing that, if $\mu_p$ denote the {\it
Segre type} of ${\Cal M}$ at $p$, the Segre sets ${\Cal S}_p^{2\mu_p}$
and ${\underline{\Cal S}}_p^{2\mu_p}$ contain an open neighborhood of
${\Cal M}$ at $p$, if $M$ is minimal at $p$ (minimality criterion due
to Baouendi, Ebenfelt and Rothschild).
 
In summary, two main steps arise in our proof, as in [BER97]
[BER99]. Step I: establishing the analyticity of jets at the level of
subsequent Segre chains and Step II: propagating analyticity up to the
maximal Segre set. In this paper, the difficulty underlying Step I is
hidden behind Artin's theorem, while Step II relies upon known
techniques of propagation. After a first reading ({\it cf.} \S11), the
reader might notice that, contrary to the very explicit iteration
process which may be endeavoured in the case of S-solvable CR maps, as
in [BER97], in the S-nondegenerate, the iteration process becomes
highly nonexplicit and requires a step by step patient induction.

\subhead \S 1.4. Separate nondegeneracy conditions \endsubhead
Now, we come to some {\it various separate assumptions on $M, f, M'$
which insure that $f$ is either S-solvable, S-finite or
S-nondegenerate.} As a main point in our definitions, the mapping $h\:
(M, p)\to_{\Cal F} (M' ,p')$ was considered as a whole object, but
some independent hypotheses on $M'$ plus other ones on $f$ are usually
made in the literature (see {\it e.g.} [BER97] [BER99] and the
references therein) and we begin by recalling some of them.

\proclaim{Proposition 1.4.1} 
{\text{\rm ([BER97])}} The formal holomorphic mapping $h\: (M,
p)\to_{\Cal F} (M' ,p')$ is $S$-solvable in each one of the following
circumstances:
\roster
\item"(i)" If $n=n'$, $M'$ is finitely nondegenerate
at $p'$ and $h$ has formal rank $n$ at $p$,
\item"(ii)" If $M\subset \C^n$, 
$M'\subset \C^{n'}$, $m\geq m'$, $M'$ is finitely nondegenerate at
$p'$ and $h$ is CR-submersive.
\endroster\endproclaim

\remark{Remarks}
1. Of course, a formal biholomorphism is CR-submersive, so (ii)
 $\Rightarrow$ (i).

\noindent
2. S-solvability of $h$ imposes furthermore a strong nondegeneracy
condition on $M'$. Indeed, one can easily see that it is necessary
that $M'$ be finitely nondegenerate at $p'$, but this is far from
being sufficient. This is why in (i) and (ii) above, $h$ is assumed to
be formally submersive on Segre varieties. Not to mention that there
exist many S-solvable mappings which are not CR submersive, see
\S13, Example 13.1.
\endremark

\proclaim{Proposition 1.4.2} 
{\text{\rm ([BER99])}} The formal holomorphic mapping $h\: (M,
p)\to_{\Cal F} (M' ,p')$ is $S$-finite in each one of the following
circumstances:
\roster
\item"(i)" If $n=n'$, $M'$ is essentially finite
at $p'$ and $h$ has formal rank $n$ at $p$,
\item"(ii)" If $n=n'$, $M'$ is essentially finite
at $p'$ and $h$ induces a finite formal map $(S_{\bar{p}}, p)
\to_{\Cal F} (S_{\bar{p}'}, p')$, 
\item"(iii)" If $M\subset \C^n$, 
$M'\subset \C^{n'}$, $m\geq m'$, $M'$ is essentially finite at $p'$
and $h$ induces a formal map $(S_{\bar{p}}, p)
\to_{\Cal F} (S_{\bar{p}'}, p')$ of generic rank equal to $m'
={\text{\rm dim}}_{\C} S_{\bar{p}'}$.
\endroster\endproclaim

A formal holomorphic map $h \: (X, p)\to_{\Cal F} (X', p')$ of complex
manifolds is said to have formal generic rank $m'={\text{\rm
dim}}_{\C} X'$ if in a local chart, a $m'\times m'$ minor of the
formal Jacobian matrix of $h$ at $p$ does not vanish identically as a
power series. The definition of finite formal maps also can be
modeled on the definition of finite holomorphic maps. Then of course
(i) $\Rightarrow$ (ii) $\Rightarrow$ (iii).

The S-finiteness of $h$ imposes a strong nondegeneracy condition on
$M'$ at $p'$: {\it for $h$ to be S-finite at $p$, it is necessary that
$M'$ be essentially finite at $p'$ but not at all sufficient} (left to
the reader). In fact, the additional conditions that are required in
Proposition 1.4.2 are all sufficient to imply that $h$ is S-finite but
they do not cover all the cases where $h$ may be S-finite, see \S13,
Example 13.2.

Let us finally remark that, although S-finite maps are not
S-nondegenerate in general, it is a fact that

\proclaim{Proposition 1.4.3} All the
formal maps $h\: (M, p)\to_{\Cal F} (M', p')$ satisfying conditions
{\rm (i)} or {\rm (ii)} of Proposition 3 or conditions {\rm (i)}, {\rm
(ii)} or {\rm (iii)} of Proposition 4 are, moreover,
S-nondegenerate. Furthermore, all the formal maps which appear in
\text{\rm[BER99]} are S-finite {\bf and} S-nondegenerate.
\endproclaim

Consequently, with our Theorem 1.2.1 (iii), we recover also all the
convergence results in [BER99] (especially Theorem 2.1 there). We can
provide here a complete independent proof of Theorem 1.2.1 (i) and
(ii) (which is almost contained in [BER97] and [BER99] respectively).

\subhead \S 1.5. Nondegeneracy conditions for generic manifolds \endsubhead
Now, we summarize the comparison between various nondegeneracy
conditions for real analytic CR manifolds. We shall say that a CR
${\Cal C}^{\omega}$ manifold is {\it S-nondegenerate at $p$} if the
identity map $i \: (M, p)\to_{\Cal F} (M,p)$ is an S-nondegenerate
formal map at $p$. Recall that $M$ is called {\it holomorphically
nondegenerate at $p$} if there does not exist a holomorphic vector
field tangent to an open piece of $M$. Also, $M$ is called {\it
finitely nondegenerate at $p$} if and only if the identity map $i \:
(M, p)\to_{\Cal F} (M,p)$ is S-solvable at $p$ and $M$ is essentially
finite at $p$ if and only if $i \: (M, p)\to_{\Cal F} (M,p)$ is
S-finite at $p$.

Assuming that $M$ is holomorphically nondegenerate, we then have:

(i) The set $\Sigma_{{\Cal F}{\Cal D}}$ of {\it finitely degenerate
points} of $M$ is a proper real analytic subvariety of $M$.

(ii) The set $\Sigma_{{\Cal N}{\Cal E}{\Cal S}{\Cal S}{\Cal F}}$ of
{\it non essentially finite points} of $M$ is a proper real analytic
subvariety of $M$.

(iii) The set $\Sigma_{{\Cal S}{\Cal D}}$ of {\it S-degenerate points}
of $M$ is a proper real analytic subvariety of $M$.

Furthermore, the following trivial inclusions
$$
\Sigma_{{\Cal F}{\Cal D}} \supset 
\Sigma_{{\Cal N}{\Cal E}{\Cal S}{\Cal S}{\Cal F}} \supset
\Sigma_{{\Cal S}{\Cal D}}
\tag 1.5.1
$$
are all strict in general (see the examples of \S13). This shows that
the condition of S-nondegeneracy is an intermediate new nondegeneracy
condition between essential finiteness and holomorphic
nondegeneracy. By the way, our technique in this paper applies well
only to S-nondegenerate maps, but we treat elsewhere the
holomorphically nondegenerate case in codimension one [MER99c].

To conclude the presentation of our main results, we would like to
mention that we can easily obtain an analog of Propositions 1.4.1 and
1.4.2 as follows.

\proclaim{Proposition 1.5.2}
The formal map $h\: (M, p)\to_{\Cal F} (M', p')$ is S-nondegenerate in
each one of the following circumstances:
\roster
\item"(i)" If $n=n'$, $M'$ is S-nondegenerate
at $p'$ and $h$ is of formal rank $n$,
\item"(ii)" If $M'$ is S-nondegenerate at $p'$ and $h$
induces a formal map $(S_{\bar{p}}, p)
\to_{\Cal F} (S_{\bar{p}'}, p')$ of generic rank equal to $m'
={\text{\rm dim}}_{\C} S_{\bar{p}'}$.
\endroster\endproclaim

Finally, it is clear that applying our Theorem 1.2.1 in all of these
situations, we thus obtain a collection of seven corollaries that we
can summarize in a

\proclaim{Theorem 1.5.3}
Under the assumptions of Propositions 1.4.1, 1.4.2 and 1.5.2, the
formal maps $h\: (M, p)\to_{\Cal F} (M', p')$ are all convergent if
$M$ is minimal at $p$.
\endproclaim

\remark{Organization of the paper}
We occupy \S2-3 with notational ingredients about Segre
chains. Despite heaviness of the underlying formalism, the geometric
picture can be easily sketched: the Segre chains happen to be just
orbits of a system of two $m$-dimensional vector fields [Su]
[Mer98]. This is why a constant use of flows of vector fields is
injected in the formalism, especially during the proofs of our main
Theorem 1.2.1, (i), (ii) and (iii). Paragraphs \S5-6 are devoted to
these proofs, using as a main tool the Theorem 1.3.2 to Artin's
theorem, a corollary which we will derive in \S12 directly from the
approximation theorem. We produce in \S13 some elementary examples to
show that many S-finite maps exist, which are not examplified by
separate assumptions on $M$, $h$, $M'$ like in Proposition 1.4.2 ({\it
cf.} [BER97] [BER99] [CPS99]) and to show that a S-finite map need not
to be S-nondegenerate in general, and vice-versa. Finally, in \S15, we
propose to the interested reader some {\it open problems}, which do
not appear explicitely in the literature, some of which are left open
by our analysis, and others of a wider class.
\endremark

\head \S2. Real analytic CR manifolds \endhead 

The next two paragraphs are devoted to a brief presentation of the
theory of Segre chains. A reader who is aware of this theory can skip
\S2 and \S3.

\subhead \S2.1. Equations \endsubhead
Let $M$ be a piece through $0$ of a ${\Cal C}^{\omega}$ {\it generic}
manifold in $\C^n$, let $m={\text{\text{\rm dim}}}_{CR} M$,
$d={\text{\rm codim}}_{\R} M$, with $m+d=n$, $\hbox{dim}_\R M=2m+n$,
and let $\rho=(\rho_1,\ldots,\rho_d)$ be a system of real analytic
defining equations for $M$ in a neighborhood $U$ of $0$ in $\C^n$,
{\it i.e.} $M=\{t\in U \:
\rho(t,\bar{t})=0\}$, $\rho(0)=0$ and $\partial\rho_1 \wedge \cdots
\wedge \partial
\rho_d(0) \neq 0$. We shall say that
$\rho$ is a {\it $d$-vectorial function}. After a sufficiently large
dilatation of the coordinates, we can assume that all the series
$\rho_j(t,\bar{t})=\sum_{\mu,\nu\in \N^n} \rho_{j,\mu, \nu} t^{\mu}
\bar{t}^{\nu}$, $\rho_{j,\mu,\nu} \in \C$, $\rho_{j, \mu,\nu}=
\bar{\rho}_{j,\nu,\mu}$, $1\leq j\leq d$,
converges uniformly in the polydisc $(4\Delta)^n$, where $\Delta$ is
the unit disc in $\C$. We let $\tau=(\bar{t})^c$ be the complexified
$\bar{t}$ variable, which is an independent variable, and we set
$\bar{\rho}(t,\tau)=\sum_{\mu,\nu\in \N^n}
\bar{\rho}_{\mu,\nu} t^{\mu} \tau^{\nu}$, so that we have
$\overline{\rho(t,\tau)}=\bar{\rho}(\bar{t},\bar{\tau})$. Also, it is
clear that there exist holomorphic coordinates
$(w,z)=(w_1,\ldots,w_m,z_1,\ldots,z_d)$ near $0\in\C^n$ vanishing at
$0\in M$ such that $T_0^c M =\C_w^m\times \{0\}$, $T_0M=\C_w^m
\times\R_x^d$, $z=x+iy$. Then $M$ is given by a system
of $d$ scalar real analytic equations $y=h(w,\bar{w},x)$ (in vectorial
notation), where $h=\sum_{\alpha,\beta,k} h_{k,\beta,\alpha} x^k
w^{\beta}
\bar{w}^{\alpha}$, $h(0)=0$, $dh(0)=0$, $h_{k,\beta,\alpha} \in \C^d$,
$\bar{h}_{k,\alpha,\beta}= h_{k,\beta,\alpha}$, $k\in \N^d$, $\beta\in
\N^m$, $\alpha\in \N^m$. Again after dilatation of $(w,\bar w,x)$, we can 
assume that the $d$-vectorial power series $h$ converges uniformly in
$(4\Delta)^{2m+n}$. Now that this choice of coordinates has been
performed, we set $\rho(t,\bar t):= y-h(w,\bar w,x)$ definitely.

In the sequel, the reference point is thought to be the origin and $p$
will denote a possibly varying point of $M$, close to the origin.

Let $\sigma$ denote the antiholomorphic involution on $\C_t^n\times
\C_{\tau}^n$ defined by $\sigma(t,\tau)=(\bar{\tau}, \bar{t})$, so
$\sigma^2={\text{\rm id}}$. To the complexification $\rho(t,\tau)$ of
$\rho(t,\bar{t})$ is canonically associated an {\it extrinsic
complexification} $M^c$ of $M$ given by $M^c=:{\Cal M}= \{(t,\tau) \in
\Delta^n \times \Delta^n \: \rho(t,\tau)=0\}
\subset \C_t^n \times \C_{\tau}^n$. Since $\partial \rho_1 \wedge
\cdots \wedge \partial \rho_d(0)\neq 0$, we have 
$\hbox{dim}_{\C} \Cal M=2m+n$. We can embed $\C^n$ in $\C_t^n \times
\C_{\tau}^n$ as the totally real plane
$\underline{\Lambda}=\{(t,\tau)\in \C^{2n} \:
\tau=\bar{t}\}$, {\it i.e.} as the graph of $t\mapsto \bar{t}$. Hence
$M$ embeds in $\underline{\Lambda}$ as $M=\{(t,\bar t) \: t\in
M\}$. Notice that $M$ also embeds in $\Cal M$, and thus $M$ can be
considered as a maximally real submanifold of ${\Cal M}$. Notice that
$\sigma(t,\bar{t})=(t,\bar{t})$, {\it i.e.} $\sigma$ fixes
$\underline{\Lambda}$ pointwisely, so $\sigma(M)= M$. If $p\in M$,
{\it i.e.} $t_p\in M\subset \C^n_t$, where $p$ is identified with its
coordinates $t_p=(t_{1,p},\ldots,t_{n,p})$, or equivalently, if $(t_p,
\bar t_p)\in M
\subset \Cal M\subset \C^n_t\times \C^n_\tau$, 
let us denote by $p^c$ the point $(t_p,\bar{t}_p)\in {\Cal M}$. Then
$p^c=\pi_t^{-1}(\{p\})\cap
\underline{\Lambda}$, where $\pi_t: \C_t^n \times \C_{\tau}^n \to 
\C_t^n$, $(t,\tau)\mapsto t$. We also put 
$\pi_{\tau}: \C_t^n \times \C_{\tau}^n
\to \C_{\tau}^n$, $(t,\tau) \mapsto \tau$, 
so that $p^c=\pi_{\tau}^{-1}(\{\bar{p}\}
\cap \underline{\Lambda})$. Using the reality of the $d$-vectorial function 
$\rho$, namely using that $\rho(t,\tau)\equiv \bar\rho(\tau,t)$, one
can easily prove that the complex manifold ${\Cal M}$ is fixed by
$\sigma$, {\it i.e.} that $\sigma({\Cal M})={\Cal M}$, and that there
is a one-to-one correspondence between germs of real analytic subsets
$\Sigma \subset M$ at $0$ and germs at $0$ of complex analytic
subvarieties $\Sigma_1 \subset {\Cal M}$ satisfying
$\sigma(\Sigma_1)=\Sigma_1$ (see [Mer98], \S2.2).
 
If we replace $y=(z-\bar{z})/2i$, $x=(z+\bar{z})/2$ in the equation
$y=h(w,\bar w,x)$ and solve this equation in terms of $z$ or of $\bar
z$, using the analytic implicit function theorem, we may obtain two
new equivalent systems of $d$-vectorial equations for $M$
$$
z=\bar{z}+i\bar{\Theta}(\bar{w},w,\bar{z}) \ \ \ \ \ {\text{\rm or}} \
\ \ \ \ \bar{z}=z-i\Theta(w,\bar{w},z),
\tag 2.1.1 
$$
with $\Theta\in \C\{w,\bar w,\bar z\}$ and thus also, two equivalent
systems of equations for ${\Cal M}$
$$
z=\xi+i\bar{\Theta}(\zeta, w,\xi) \ \ \ \ \ {\text{\rm or}} \ \ \ \ \
\xi=z-i\Theta(w,\zeta,z),
\tag 2.1.2 
$$
where $\zeta=(\bar{w})^c$ and $\xi=(\bar{z})^c$. Since $h(0)=0$ and
$dh(0)=0$, then we have also $\bar\Theta(0)=0$ and $d\bar
\Theta(0)=0$. After a new dilatation of the coordinates, we can assume
the convergence in $(4\Delta)^{2m+n}$ of the $d$-vectorial function
$\bar{\Theta}$ above. The fact that these two systems of equations
define the same manifold ${\Cal M}$ (or, equivalently, that
$\sigma({\Cal M})={\Cal M}$) is reflected by the following relations
that are obtained by replacing one system of equations of ${\Cal M}$
into the other
$$
\Theta(w,\zeta,z)\equiv \bar{\Theta}
(\zeta,w,z-i\Theta(w,\zeta,z)) \ \ {\text{\rm and}} \ \
\bar{\Theta}(\zeta,w,\xi) \equiv \Theta(w,\zeta,\xi
+i\bar{\Theta}(\zeta,w,\xi)).
\tag 2.1.3
$$

\subhead \S 2.2. Vector fields \endsubhead
We say that $L$ is an $m$-vector field if $L$ is given by the span of
a set of $m$ independent vector fields $L_1,\ldots,L_m$ over, say,
$\Delta^n \subset \C^n$ or $\Delta^{2n} \subset \C^{2n}$, {\it such
that $L_1,\ldots,L_m$ commute with each other}.

For instance, the complexification ${\Cal M}=M^c$ admits complexified
(1,0) and (0,1) tangent $m$-vector fields which can be written
explicitely as
$$
{\Cal L}=\frac{\partial }{\partial w}+i\bar{\Theta}_w (\zeta,w,\xi)
\frac{\partial }{\partial z} \ \ {\text{\rm and}} \ \
\underline{\Cal L} =\frac{\partial }{\partial \zeta}-i\Theta_{\zeta}
(w,\zeta,z) \frac{\partial }{\partial \xi}
\tag 2.2.1
$$
in vectorial notation: here $\Cal L=(\Cal L_1,\ldots, \Cal L_m)$ and
$\underline{\Cal L}= (\underline{\Cal L}_1,\ldots,\underline{\Cal
L}_m)$ each form a system of $m$-vectorial vector fields over
$\Delta^{2n}$. These two $m$-vector fields form together a couple
$\SS= \{{\Cal L},
\underline{\Cal L}\}$ of $m$-vector fields that are clearly the
complexifications of the canonical representatives for a basis of
$T^{1,0}M$ and $T^{0,1} M$ respectively, namely
$$
L=\frac{\partial }{\partial w}+i\bar{\Theta}_w(\bar{w},w,\bar{z})
\frac{\partial }{\partial z} \ \ {\text{\rm and}} \ \
\bar{L} = \frac{\partial }{\partial \bar{w}}-i\Theta_{\bar{w}} (w,\bar{w},z)
\frac{\partial }{\partial \bar{z}},
\tag 2.2.2
$$
that is to say $L^c={\Cal L}$ and $\bar{L}^c=\underline{\Cal L}$.

Let us recall the following fundamental facts (see [Mer98], \S2). The
system $\{L,\bar{L} \}$ is orbit-minimal on $M$ if and only if the
system $\{{\Cal L}, \underline{\Cal L}\}$ is orbit-minimal on ${\Cal
M}$ and, more generally, ${\Cal O}_{L, \bar{L}}(M,p)^c= {\Cal
O}_{{\Cal L},\underline{\Cal L}}({\Cal M}, p^c)$. The main fact is
that the family of complexified Segre and conjugate Segre varieties
${\Cal S}_{\tau_p}$ and $\underline{\Cal S}_{t_p}$ form families of
integral manifolds for ${\Cal L}$ and $\underline{\Cal L}$
respectively which induce their canonical flow foliation. Indeed, this
can be observed after writing explicitely the complexifications of the
Segre and of the conjugate Segre varieties, namely:
$$
S_{\bar{t}_p} : z=\bar{z}_p + i\bar{\Theta}(\bar{w}_p,w,\bar{z}_p) \ \
{\text{\rm and}} \ \ \overline{S}_{t_p} :
\bar{z}=z_p-i\Theta(w_p,\bar{w},z_p)
\tag 2.2.3
$$
whose complexifications can be seen and written as follows:
$$
\aligned
& {\Cal S}_{\zeta_p,\xi_p}:
\zeta=\zeta_p, \xi=\xi_p, z=\xi_p+i\bar{\Theta}(\zeta_p,
w,\xi_p) \ \ {\text{\rm and}} \\ &
\underline{{\Cal S}}_{w_p,z_p}: 
w=w_p, z=z_p, \xi=z_p-i\Theta(w_p,\zeta,z_p),
\endaligned
\tag 2.2.4
$$
and thus it is clear from (2.2.4) that ${\Cal L} {\Cal
S}_{\zeta_p,\xi_p}\equiv 0$ and $\underline{\Cal L} \underline{\Cal
S}_{w_p,z_p}\equiv 0$. For dimensional reasons ($\hbox{dim}_\C {\Cal
S}_{\zeta_p,\xi_p}= \hbox{dim}_\C \underline{{\Cal S}}_{w_p,z_p}=m$),
these complexified Segre varieties coincide with the leaves of the
flow foliations induced by $\Cal L$ and by $\underline{\Cal L}$ on
$\Cal M$, respectively.

\head \S3. Segre chains\endhead
 
\subhead \S3.1. Definitions \endsubhead
The {\it orbit $k$-chains} of the pair $\{{\Cal L}, \underline{\Cal
L}\}$ of $m$-vector fields on ${\Cal M}$ are called {\it Segre
$k$-chains}. The orbit $k$-chains are almost implicitely defined by
Sussmann in [Su] and the reader is referred to [Mer98] for the general
construction in the analytic category. Here, we summarize how they can
be constructed.

First, let us write ${\Cal L}=({\Cal L}_1,\ldots,{\Cal L}_m)$ and let
$w\in \C^m$. We have already observed that the ${\Cal L}_j$'s
commute. Consequently, if we consider the multiple flow mapping
$$
\C^m\times
\Delta^{2n} \ni (w_1,\ldots,w_m,p)\mapsto
\exp (w_m {\Cal L}_m) \circ \cdots \circ \exp (w_1 {\Cal L}_1)(p)\in 
\Delta^{2n},
\tag 3.1.1
$$
which is defined over a certain domain of $\C^m\times
\Delta^{2n}$, we shall have immediately:
$$
\exp(w_{\varpi(m)} {\Cal L}_{\varpi(m)})\circ \cdots \circ
\exp(w_{\varpi(1)} {\Cal L}_{\varpi(1)})(p)=
\exp (w_m {\Cal L}_m) \circ \cdots \circ \exp(w_1 {\Cal L}_1)(p),
\tag 3.1.2
$$ 
for every permutation $\varpi: \{1,2,\ldots,m\} \to \{1,2,\ldots,m\}$.
The $m$-vector field $\underline{\Cal L}$ satisfies the same
property. We shall simply denote these two multiple flow maps by
$(w,p)\mapsto {\Cal L}_w(p)$ and by $(\zeta,p)\mapsto \underline{\Cal
L}_{\zeta}(p)$. Thus, we can formally work with multiple vector fields
as if they were usual vector fields, {\it i.e.} as if $m=1$.

Now, let us denote shortly $\SS=\{{\Cal L}, \underline{\Cal
L}\}$. Then ${\text{\rm rk}} \ \SS=2m$ at each point of $\Delta^{2n}
\cap {\Cal M}$. The {\it codimension} of $\SS$ in ${\Cal M}$ equals
${\text{\rm dim}}_{\C} {\Cal M}-2m=d={\text{\rm codim}}_{\R} M$. We
fix a {\it reference} point $p\in M$, $p^c\in {\Cal M}$, close to the
origin 0 (which is the {\it central} point in our coordinates
$(t,\tau)$). It will be interesting to let this point $p$ vary in
order to consider the various orbits of different points close to the
origin. In this paragraph, the reference point should not be confused
with the origin.

If $k\in \N_*$, $\L^k=(L^1,\ldots,L^k)\in \SS^k$,
$w_{(k)}=(w_1,\ldots,w_k)\in \C^k$, $p\in (\frac{1}{2}\Delta)^{2n}\cap
{\Cal M}$, let us denote $\L^k_{w_{(k)}}(p)= L_{w_k}^k \circ \cdots
\circ L_{w_1}^1(p)$ whenever the composition is defined. By the way, after
bounding $k\leq 3(2n)$, it is clear that $\L_{w_{(k)}}^k(p)\in
\Delta^{2n} \cap {\Cal M}$ whenever $w_{(k)}\in (\delta \Delta^m)^k$,
$p\in (\frac{1}{2} \Delta)^{2n} \cap {\Cal M}$, $k\leq 3(2n)$, for
$\delta >0$ small enough.

Because ${\text{\rm Card}} \ \SS=2$, only two different $\SS$-chains
exist. They can be naturally called Segre $k$-chains (we choose this
denomination instead of ``Segre sets'', because they come together
with compositions of flow paps, which are concatenations of
holomorphic maps). These two families of Segre $k$-chains of $p=
(t_p,\tau_p)\in {\Cal M}$ are defined by
$$
\aligned
& {\Cal S}_{\tau_p}^{2j}=\{\underline{\Cal L}_{\zeta_j}\circ {\Cal
L}_{w_j}\circ\cdots\circ\underline{\Cal L}_{\zeta_1}\circ{\Cal
L}_{w_1}(p)\: w_1,\zeta_1,\ldots,w_j,\zeta_j\in \delta\Delta^m\}\\ &
{\Cal S}_{\tau_p}^{2j+1}=\{{\Cal L}_{w_{j+1}}\circ \underline{\Cal
L}_{\zeta_j}
\circ\cdots\circ\underline{\Cal L}_{\zeta_1}\circ{\Cal
L}_{w_1}(p)\: w_1,\zeta_1,\ldots,\zeta_j,w_{j+1}\in
\delta\Delta^m\}\\
&
\underline{\Cal S}_{t_p}^{2j}=\{{\Cal
L}_{w_j}\circ \underline{\Cal L}_{\zeta_j}\circ\cdots\circ {\Cal
L}_{w_1}\circ \underline{\Cal L}_{\zeta_1}(p) \:
\zeta_1,w_1,\ldots,\zeta_j,w_j\in \delta\Delta^m\}\\
& {\Cal S}_{t_p}^{2j+1}=\{\underline{\Cal L}_{\zeta_{j+1}}\circ {\Cal
L}_{w_j}\circ\cdots\circ {\Cal L}_{w_1}\circ \underline{\Cal
L}_{\zeta_1}(p) \:
\zeta_1,w_1,\ldots,w_j,\zeta_{j+1}\in \delta\Delta^m\},
\endaligned
\tag 3.1.3
$$
for $k=2j$ or $k=2j+1$, $j\in \N$, $k\leq 3(2n)$. Of course, ${\Cal
S}_{\tau_p}^k \subset {\Cal M}$ and $\underline{\Cal S}_{t_p} \subset
{\Cal M}$.

The actions of the vectorial $m$-flows of ${\Cal L}$ and of
$\underline{\Cal L}$ on a point $p\in \Cal M$ with coordinates
$(w_p,z_p,\zeta_p,\xi_p)\in \Cal M$ are simply given by
$$
\aligned
{\Cal L}_w(w_p,z_p=\xi_p+ & i\bar{\Theta}(\zeta_p,w_p, \zeta_p),
\zeta_p,\xi_p)= 
\\
& =(w_p+w,z=\xi_p+i\bar{\Theta}(\zeta_p, w_p+w, \zeta_p),
\zeta_p, \xi_p)\\
\underline{\Cal L}_{\zeta}(w_p, z_p, \zeta_p,
\xi_p= 
& z_p-i\Theta(w_p, \zeta_p,z_p))=
\\
& =(w_p,z_p, \zeta+\zeta_p,
\xi=z_p-i\Theta(w_p,\zeta+\zeta_p, z_p)).
\endaligned
\tag 3.1.4
$$

We recall the property $\sigma({\Cal L}_w(q))=\underline{\Cal
L}_{\bar{w}}(\sigma(q))$, which can be seen easily after applying
$\sigma$ to the members of eq. (3.1.4), so that more generally
$\sigma({\Cal S}_{\tau_p}^k)=\underline{\Cal S}_{\bar{\tau}_p}^k$ and
$\sigma(\underline{\Cal S}_{t_p}^k)={\Cal S}_{\bar{t}_p}^k$. This
property easily implies that the {\it a priori} different two
minimality types and multitypes of $\SS$ must coincide. Accordingly,
let us recall how the Segre type and multitype of ${\Cal M}$ at $0$
are defined.

First, let us point out that, if we are given $f: X\to Y$ a
holomorphic map of connected complex manifold, then there exists a
proper complex subvariety $Z\subset X$ with ${\text{\rm dim}}_{\C} Z <
{\text{\rm dim}}_{\C} X$ and an integer $r$ such that ${\text{\rm
rk}}_{\C,p}(f)=r=\max_{q\in X} {\text{\rm rk}}_{\C,q} (f)$ for all
$p\in X\backslash Z$. We shall denote this integer, the {\it generic
rank} of $f$, by ${\text{\rm genrk}}_{\C}(f)$.

As ${\text{\rm Card}}\ \SS=2$, any alternating $k$-tuple of
$m$-vectors $\L^k\in \SS^k$ can be written uniquely $(\cdots, {\Cal
L}, \underline{\Cal L}, {\Cal L})$ or $(\cdots,
\underline{\Cal L}, {\Cal L}, \underline{\Cal L})$. Accordingly, we shall write
${\Cal L} \underline{\Cal L}^k$ to denote $({\Cal L}\underline{\Cal
L})^j$ if $k=2j$ is even and $\underline{\Cal L} ({\Cal L}
\underline{\Cal L})^j$ if $k=2j+1$ is odd. Of course, we introduce also
the similar notation for $\underline{\Cal L} {\Cal L}^k$. Also, we
denote by $\Gamma_{\underline{\Cal L} {\Cal L}^k}$ the map $(\delta
\Delta^m)^k\to \Delta^{2n}$, $w_{(k)} \mapsto \underline{\Cal L}{
\Cal L}_{w_{(k)}}^k(0)$, where $w_{(k)}=(w_1,\ldots,w_k)$.

In these notations, ${\Cal S}_{\tau_p}^k=\Gamma_{\underline{
\Cal L}{\Cal L}^k}((\delta
\Delta^m)^k)$ and $\underline{\Cal S}_{t_p}^k=\Gamma_{{\Cal L}
\underline{\Cal L}^k}((\delta \Delta^m)^k)$. By slight abuse of notation
we shall also denote $\Gamma_{\underline{\Cal L} {\Cal L}^k}$,
$\Gamma_{{\Cal L}\underline{\Cal L}^k}$ by $\Gamma_{\underline{\Cal L}
{\Cal L}}^k$, $\Gamma_{{\Cal L}
\underline{\Cal L}}^k$. Whenever we will be working with 
coordinates, the central point will be for us the origin and $p$ will
denote a varying point of ${\Cal M}$. When we state an invariant
theorem, $p$ is the reference point that we pick in $\Cal M$ or in
$\Cal M$ near the origin.

The {\it Segre multitype of ${\Cal M}$ at $p$} is defined to be the
$\mu_p$-tuple $(m,m,e_1,\ldots,e_{\kappa_p})$, $\mu_p=2+\kappa_p$,
being the {\it Segre type of} ${\Cal M}$ at $p$ satisfying

1. ${\text{\rm genrk}}_{\C} (\Gamma_{\underline{\Cal L} {\Cal L}}^k)
=2m+e_1+\cdots+e_k=2m+e_{\{k\}}$, $2\leq k\leq \kappa_p$ and

2. ${\text{\rm genrk}}_{\C} (\Gamma_{\underline{\Cal L} {\Cal
L}}^{\mu_p})= {\text{\rm genrk}}_{\C} (\Gamma_{\underline{\Cal L}
{\Cal L}}^{\mu_p+1})= 2m+ e_{\{\kappa_p\}}$, $\mu_p=2+\kappa_p$.

\subhead \S3.2. Properties \endsubhead
All these integers are uniquely and invariantly defined (because the
two canonical Segre foliations are biholomorphically invariant). Of
course, ${\text{\rm genrk}}_{\C} (\Gamma_{\underline{\Cal L} {\Cal
L}}^1)=m$ and ${\text{\rm genrk}}_{\C} (\Gamma_{\underline{\Cal L}
{\Cal L}}^2)=2m$. Using elementary properties of the flow maps, one
can show easily ({\it cf.} [Mer98]) that ${\text{\rm genrk}}_{\C}
(\Gamma_{\underline{\Cal L} {\Cal L}}^{\mu_p+k})= 2m+e_{\{\kappa_p\}}$
for all $k\geq 1$ and that $e_1>0,\ldots,e_{\kappa_p}>0$. The integer
$\mu_p$ satisfying properties 1 and 2 above is called the {\it Segre
type} of ${\Cal M}$ at $p$ and the $\mu_p$-tuple
$(m,m,e_1,\ldots,e_{\kappa_p})$ is called the {\it Segre multitype} of
${\Cal M}$ at $p$. There is a similar definition of multitype for the
maps $\Gamma_{{\Cal L}\underline{\Cal L}}^k$ instead of
$\Gamma_{\underline{\Cal L}{\Cal L}}^k$, but this definition yields
the same integers, because of the property $\sigma({\Cal
L}\underline{\Cal L}_{w_{(k)}}^k(q))=
\underline{\Cal L} {\Cal L}_{\bar{w}_{(k)}}^k(\sigma(q))$.

Let us summarize the properties of Segre chains in ${\Cal M}$. For the
definition and the presentation of the concept of {\it orbit} of a
system of vector fields, see [Su] [Mer98].

\proclaim{Theorem 3.2.1}
Let ${\Cal M}=M^c$, $m={\text{\rm dim}}_{CR} M$, $d={\text{\rm
codim}}_{\R} M \geq 1$, ${\text{\rm dim}}_{\R} M= 2m+d={\text{\rm
dim}}_{\C} {\Cal M}$ and let $p^c\in {\Cal M}$, $p^c=
\pi_t^{-1}(\{p\}) \cap \underline{\Lambda}$. Then there
exist an invariant integer $\mu_p\in \N_*$, $3\leq \mu_p \leq d+2$,
the {\rm Segre type} of ${\Cal M}$ at $p^c$ and
$(m,m,e_1,\ldots,e_{\kappa_p})$, the {\rm Segre multitype} of ${\Cal
M}$ at $p^c$, $\kappa_p=\mu_p-2$, and $w_{(\mu_p)}^*\in
(\delta\Delta^m)^{\mu_p}$ with $w_{\mu_p}^*=0$ and a neighborhood
${\Cal W}^*$ of $w_{(\mu_p)}^*$ in $(\delta\Delta^m)^{\mu_p}$ such
that, putting $\underline{w}_{(\mu_p-1)}^*=
(-w_{\mu_p-1}^*,\ldots,-w_1^*)$, we have:

$1)$ \ $\Gamma_{\underline{\Cal L}{\Cal L}}^{2\mu_p-1}({\Cal
W}^*\times\{\underline{w}_{(\mu_p-1)}^*\})={\Cal O}_p^c=$ constitutes
a piece ${\Cal O}_p^c$ of the complexified CR orbit of $M$ through
$p^c=$ ${\Cal N}=$ a piece ${\Cal N}$ of the orbit ${\Cal O}_{{\Cal
L}, \underline{\Cal L}} ({\Cal M}, p^c)$ through $p^c${\rm ;}

$2)$ \ $2m+e_1+\cdots+e_{\kappa_p}=2m+e_{\{\kappa_p\}}={\text{\rm
dim}}_{\C} {\Cal O}^c_p={\text{\rm dim}}_{\C} {\Cal O}_{{\Cal L},
\underline{\Cal L}} ({\Cal M}, p^c)= {\text{\rm dim}}_{\R} {\Cal
O}_p${\rm ;}

$3)$ \ ${\text{\rm genrk}}_{\C} (\Gamma_{\underline{\Cal L}{\Cal
L}}^{k+2})= 2m+e_{\{k\}}=2m+e_1+\cdots+ e_k= {\text{\rm genrk}}_{\C}
(\Gamma_{{\Cal L}\underline{\Cal L}}^{k+2})$, $0\leq k \leq
\kappa_p${\rm ;}

$4)$ \ ${\text{\rm genrk}}_{\C} (\Gamma_{\underline{\Cal L}{\Cal
L}}^{k+2})=2m+ e_{\{\kappa_p\}} ={\text{\rm genrk}}_{\C}
(\Gamma_{{\Cal L}\underline{{\Cal L}}}^{k+2})$, $\kappa_p\leq k\leq
3(2n)-2$.
\endproclaim
 
It is not in this form that we shall use the main properties of Segre
chains. A more appropriate theorem which is suitable for our purposes
can be stated as follows. This theorem will be mainly used at the very
end of the proof of our main Theorem 1.2.1, precisely in Assertion
8.4. More effectively, we will use this theorem in the case of a
minimal generic manifold, {\it i.e.} we will use its Corollary 3.2.6
below.

\proclaim{Theorem 3.2.2} 
If $\mu_p$ denotes the Segre type of ${\Cal M}$ at $p^c$, then there
exist a $\mu_p$-tuple $w_{(2\mu_p)}^*\in (\delta\Delta^m)^{2\mu_p}$
and a neighborhood ${\Cal W}^*$ of $w_{(2\mu_p)}^*$ in
$(\delta\Delta^m)^{2\mu_p}$ such that that $\Gamma_{\underline{\Cal L}
{\Cal L}}^{2\mu_p} (w_{(2\mu_p)}^*)$ $=p$ and, moreover
$$
{\text{\rm rk}}_{\C,w_{(2\mu_p)}^*}(\Gamma_{
\underline{\Cal L}{\Cal L}}^{2\mu_p})=
{\text{\rm dim}}_{\C} {\Cal O}_{{\Cal L}, \underline{\Cal L}}({\Cal
M}, p^c),
\tag 3.2.3
$$
$$
{\text{\rm rk}}_{\C, w_{(2\mu_p)}^*} (\pi_{\tau} \circ
\Gamma_{\underline{\Cal L}{\Cal L}}^{2\mu_p})=
{\text{\rm rk}}_{\C, w_{(2\mu_p)}^*} (\pi_t \circ
\Gamma_{\underline{\Cal L}{\Cal L}}^{2\mu_p})=
{\text{\rm dim}}_{\C} {\Cal O}_p^{i_c}= m+e_{\{\kappa_p\}},
\tag 3.2.4
$$
where ${\Cal O}_p^{i_c}\subset \C_t^n$ is the intrinsic
complexification of a piece ${\Cal O}_p$ of ${\Cal O}_{CR}(M,p)$.
\endproclaim

\remark{Remark} 
The rank properties (3.2.3) and (3.2.4) above follow in fact more or
less directly from the geometric construction in Sussmann's theory of
orbits. Furthermore, Theorem 3.2.1 holds for $\delta >0$ arbitrarily
small.
\endremark\medskip

In particular, if $M$ is generic and {\it minimal} at $p$, we have
$$
{\text{\rm rk}}_{\C, w_{(2\mu_p)}^*} (\pi_{\tau} \circ
\Gamma_{\underline{\Cal L}{\Cal L}}^{2\mu_p})=
{\text{\rm rk}}_{\C, w_{(2\mu_p)}^*}(\pi_t \circ
\Gamma_{\underline{\Cal L}{\Cal L}}^{2\mu_p})=n.
\tag 3.2.5
$$
Hence the appropriate characterization of minimality for our purposes
in this article:
 
\proclaim{Corollary 3.2.6} 
If $M$ is minimal at $p$, for each $\delta >0$, there exists a
$\mu_p$-tuple $(\zeta_{\mu_p}^*,
w_{\mu_p}^*,\ldots,\zeta_1^*,w_1^*)\in (\delta\Delta^m)^{2\mu_p}$ such
that $\underline{\Cal L}_{\zeta_{\mu_p}^*} \circ {\Cal
L}_{w_{\mu_p}^*}\circ \cdots \circ \underline{\Cal L}_{\zeta_1^*}
\circ {\Cal L}_{w_1^*}(p^c)=p^c$ and such that the
ranks of the two mappings
$$
(\zeta_{\mu_p},w_{\mu_p},\ldots,\zeta_1,w_1) \mapsto
\pi_t \ {\text{\rm or}} \ \pi_{\tau}
(\underline{\Cal L}_{\zeta_{\mu_p}} \circ {\Cal L}_{w_{\mu_p}}
\circ \cdots \circ \underline{\Cal L}_{\zeta_1} \circ {\Cal L}_{w_1}(p^c))
\tag 3.2.7
$$
at the point $(\zeta_{\mu_p}^*, w_{\mu_p}^*,\ldots,\zeta_1^*,w_1^*)$
are both equal to $n$.
\endproclaim

\remark{Remark} Also, there exists a $\mu_p$-tuple
$(w_{\mu_p}^{**}, \zeta_{\mu_p}^{**},\ldots,w_1^{**},\zeta_1^{**})\in
(\delta\Delta^m)^{2\mu_p}$ such that $\underline{\Cal
L}_{\zeta_{\mu_p}^{**}} \circ {\Cal L}_{w_{\mu_p}^{**}}
\circ \cdots \circ \underline{\Cal L}_{\zeta_1^{**}} 
\circ {\Cal L}_{w_1^{**}}(p^c)=p^c$
and such that the ranks of the two mappings
$$
(w_{\mu_p}, \zeta_{\mu_p},\ldots,w_1,\zeta_1) \mapsto
\pi_{\tau} \ {\text{\rm or}} \ \pi_t ({\Cal L}_{w_{\mu_p}} \circ 
\underline{{\Cal L}}_{\zeta_{\mu_p}} \circ \cdots \circ
{\Cal L}_{w_1} \circ \underline{\Cal L}_{\zeta_1} (p^c))
\tag 3.2.8
$$
are both equal to $n$ at the point $(w_{\mu_p}^{**},
\zeta_{\mu_p}^{**},\ldots,w_1^{**},\zeta_1^{**})$. In fact, thanks to
the property $\sigma({\Cal L}_w(p^c))=\underline{\Cal
L}_{\bar{w}}(\sigma(p^c))$, it easy to see that one can choose
$$
(w_{\mu_p}^{**}, \zeta_{\mu_p}^{**},\ldots,w_1^{**},\zeta_1^{**}):=
(\bar{\zeta}_{\mu_p}^*,
\bar{w}_{\mu_p}^*,\ldots,\bar{\zeta}_1^*,\bar{w}_1^*).
\tag 3.2.9
$$
\endremark\medskip
 
This completes the presentation of Segre chains.

\head \S4. The mapping \endhead 

Let $M'$ be a second real analytic CR manifold, which is generic in
$\C^{n'}$, with ${\text{\rm dim}}_{CR}M'=m'\geq 1$ and ${\text{\rm
codim}}_{\R} M'=d'\geq 1$, $m'+d'=n'$, $\dim_\R M'=2m'+d'$, and let
$p'\in M'$. As for $M$, there exist holomorphic coordinates $(w',z')$
vanishing at $p'$ in which $M'$ is given by two equivalent systems of
$d'$ scalar equations:
$$
z'=\bar{z}'+i\bar{\Theta}'(\bar{w}',w',\bar{z}') \ \ \ \ \ {\text{\rm
or}}
\ \ \ \ \ \bar{z}'=z'-i\Theta'(w',\bar{w}',z'),
\tag 4.1
$$
with $\Theta'$ converging in $(4\Delta)^{2m'+d'}$ and also, two
equivalent systems of $d'$ scalar equations for its extrinsic
complexification ${\Cal M}'$:
$$
z'= \xi'+i\bar{\Theta}'(\zeta',w',\xi') \ \ \ \ \ {\text{\rm or}} \ \
\ \ \
\xi'=z'-i\Theta'(w',\zeta',z').
\tag 4.2
$$
Let now $h:(\C^n,p)\mapsto_{\Cal F} (\C^{n'},p')$,
$h(t)=(\sum_{\alpha\in\N_*^n} h_{1,\alpha} t^{\alpha},\ldots,
\sum_{\alpha\in \N_*^n} h_{n',\alpha} t^{\alpha})$ 
be a formal holomorphic mapping, in our coordinates in which $p=0$,
$p'=0$, $h(0)=0$. By definition, $h$ is called {\it invertible at} $0$
(as in (i) of Proposition 1.5.2), if $n'=n$ and if we have ${\text{\rm
det}} (h_{i,\1_j})_{1\leq i,j\leq n}\neq 0$, if $\1_j$ denotes the
$n$-tuple $(0,\ldots,0,1,0,\ldots,0)$ with $1$ at the place of the
$i$-th digit.

We write $h=(h_1,\ldots,h_{n'})=(g_1,\ldots,g_{m'},f_1,\ldots,f_{d'})
=(g,f)$ according to the splitting of coordinates $(w',z')$ which is
coherent with the disposition of the complex tangent space to $M'$,
$T_0^cM'=\C_{w'}^{m'}\times \{0\}$. In order to provide an aesthetical
presentation of all the formalism used throughout the article, we
shall maintain a uniform {\it vectorial convention}, most often
omitting certain superfluous indices. For instance, we shall write
${\text{\rm det}} (\underline{\Cal L} \bar{g})$ instead of ${\text{\rm
det}} (\underline{\Cal L}_j
\bar{g}_k)_{1\leq j,k\leq m}$ (if $n=n'$, $m=m'$).

Without loss of generality, we can assume that the coordinates $(w,z)$
for $M$ and $(w',z')$ for $M'$ are {\it normal}, {\it i.e.} that
$\Theta(0,\zeta,z)\equiv 0$, $\Theta(w,0,z)\equiv 0$ and that
$\Theta'(0,\zeta',z')\equiv 0$, $\Theta'(w',0,z')\equiv 0$.

Our main assumption is that $h: (M,0) \to_{{\Cal F}} (M',0)$ maps
$(M,0)$ into $(M',0)$ {\it formally}. This means that the following
two equivalent formal identities (between formal power series) hold
$$
\aligned
& f(w,z)\equiv \left[
\bar{f}(\bar{w},\bar{z})+i\bar{\Theta}'(\bar{g}(\bar{w},\bar{z}), 
g(w,z),\bar{f}(\bar{w}, \bar{z}))
\right]_{\bar{z}:=z-i\Theta(w,\bar{w},z)}\\
&
\bar{f}(\bar{w},\bar{z})\equiv
\left[f(w,z)-i\Theta'(g(w,z),\bar{g}(\bar{w},\bar{z}),f(w,z))\right]_{
z:=\bar{z}+i\bar{\Theta}(\bar{w},w,\bar{z})}
\endaligned
\tag 4.3
$$
in $\C\dl w,\bar{w},z\dr$ and in $\C\dl \bar w,w,\bar{z}\dr$
respectively. Consequently, $h$ induces a formal map
$h^c:=(h,\bar{h})$, $h^c(t,\tau)=(h(t),\bar{h}(\tau))=(g(w,z),f(w,z),
\bar{g}(\zeta,\xi), \bar{f}(\zeta,\xi))$ between $({\Cal M},0)$ and
$({\Cal M}',0)$, which can be constructed just by complexifying the
two formal identities in (4.3). We obtain thus two equivalent formal
identities (between formal power series):
$$
\aligned
& f(w,z)=\left[\bar{f}(\zeta,\xi)+i\bar{\Theta}'(\bar{g}(\zeta,\xi),
g(w,z),\bar{f}(\zeta,\xi))\right]_{\xi:=z-i\Theta(w,\zeta,z)}\\ &
\bar{f}(\zeta,\xi)=\left[f(w,z)-i\Theta'(g(w,z),\bar{g}(\zeta,\xi),
f(w,z))\right]_{z:=\xi+i\bar{\Theta}(\zeta,w,\xi)}
\endaligned
\tag 4.4
$$
in $\C\dl w, \zeta, \xi\dr$ and in $\C\dl \zeta, w,\xi\dr$
respectively, after replacing $\xi$ by $z-i\Theta(w,\zeta,z)$ in the
first equation of (4.4) and $z$ by $\xi+i\bar{\Theta}(\zeta,w,\xi)$ in
the second equation. Then (4.4) means that $h^c({\Cal M}, 0)
\subset_{\Cal F} ({\Cal M}',0)$, {\it i.e.} that $h^c$ maps $(\Cal M,
0)$ formally into $(\Cal M', 0)$.

In principle, {\it both the triples $(\zeta,w,z)$ and $(w,\zeta,\xi)$}
can be chosen as coordinates on ${\Cal M}$. {\it There is no canonical
or preferred choice}. Consequently, we will opt here for a systematic
twofold presentation of everything, in coherence with the twofold
theory of Segre chains which was built [Mer98]. Furthermore, we shall
henceforth work only with the complexified map $h^c: ({\Cal M}, 0)
\to_{\Cal F} ({\Cal M}',0)$. Let us finish with a
 
\proclaim{Lemma 4.5} 
Let $h^c=(h,\bar{h})\: (\Cal M,0)\to_{\Cal F} (M',0)$ as above,
$h^c(t,\tau)=(h(t),\bar{h}(\tau))$. Then $\sigma'\circ h^c=h^c \circ
\sigma$. Also, notice that $h \circ \underline{\Cal L}_{\zeta} 
(q(x)) \equiv h(q(x))$ in $\C\dl x,\zeta\dr$ and that $\bar{h} \circ
{\Cal L}_w (q(x)) \equiv \bar{h}(q(x))$ in $\C\dl x,w\dr$ for any
formal series $q(x)\in \C\dl x\dr^{2n}$, $q(0)=0$.
\endproclaim

\demo{Proof} The first property is trivial:
$$
\sigma'\circ h^c(t,\tau)=\sigma'(h(t),\bar h(\tau))=
(h(\bar \tau), \bar h(\bar t))= h^c(\bar \tau, \bar t)= h\circ \sigma
(t,\tau).
\tag 4.6
$$
The second property follows easily from eq. (3.1.4), from which we see
that $\pi_t\circ \underline{\Cal L}_\zeta (q(x))=\pi_t(q(x))$ and
$\pi_\tau \circ \Cal L_w (q(x))=\pi_\tau(q(x))$ and from $\bar
h(r(w,x))=\bar h \circ \pi_\tau (r(w,x))$ and $h(r(\zeta,x))=h\circ
\pi_t(s(\zeta,x))$ for any two formal power series $r(w,x) \in \C\dl
w,x \dr$ and $s(\zeta,x)\in \C\dl \zeta, x\dr$.
\qed\enddemo

\head \S5. Differentiations \endhead

The purpose of this paragraph is to prove Proposition 1.5.2 (i),
namely to establish that any formal invertible holomorphic mapping
between S-nondegenerate CR manifolds is S-nondegenerate (Proposition
5.13 in this paragraph). We thus work in the equidimensional case,
{\it i.e.} with $m'=m$, $d'=d$ and $n'=n$.

We first derive a convient expression of the {\it fundamental family
of identities} written in eq. (1.1.2) as follows.

Let us consider the derivations $\underline{\Cal
L}^\beta=\underline{\Cal L}_1^{\beta_1} \cdots
\underline{\Cal L}_{n_1}^{\beta_{n-1}}$. Then applying all these derivations of
any order to the identity $\bar \rho'(\bar h(\tau),h(t))$, {\it i.e.}
to
$$
\bar{f}(\zeta, \xi) \equiv f(w,z) -i\sum_{\gamma\in \N_*} 
\bar g(\zeta, \xi)^\gamma \ \Theta_\gamma'(g(w,z), f(w,z)),
\tag 5.2
$$
as $(w,z,\zeta,\xi)\in \Cal M$, it is well-known that we obtain an
infinite family of formal identities that we collect in an independent
statement that we will reprove quickly below, for completeness ({\it
cf.} [BR88] [BR90] [BER97] [BERbk] [CPS2] [Mer99c]). Let
$\Theta_{{\zeta'}^{\beta}}'(w',\zeta',z')$ denote
$\partial_{\zeta'}^{\beta} \Theta'(w',\zeta',z')$, $\beta\in \N^m$,
$\partial_{\zeta'}^{\beta}=\partial_{\zeta_1'}^{\beta_1} \cdots
\partial_{\zeta_m'}^{\beta_m}$.

\proclaim{Proposition 5.1}
Let $h: (M,0) \to_{{\Cal F}} (M',0)$ be a formal biholomorphism
between CR generic ${\Cal C}^{\omega}$ manifolds in $\C^n$. For every
$\beta\in \N^m$, there exists a collection of $d\times d$ matrices of
universal polynomial $\underline{u}_{\beta,\gamma}$, $|\gamma|\leq
|\beta|$ in $m N_{m,|\beta|}$ variables, where $N_{m,|\beta|}={(\vert
\beta \vert+ m)!\over \vert \beta \vert ! \ m!}$, and holomorphic
$\C^d$-valued functions $\underline{\Omega}_{\beta}$ in
$(2n-d+n'N_{n,|\beta|})$ variables near $0\times 0\times 0 \times
(\partial_{\xi}^{\alpha^1}\partial_{\zeta}^{
\gamma^1}\bar{h}(0))_{|\alpha^1|+ |\gamma^1| \leq
|\beta|}$ in $\C^m\times
\C^m \times \C^d \times \C^{n'N_{n,|\beta|}}$ such that
$$
\aligned
\underline{\theta}_{\beta}'(w,\zeta,\xi)
& :=\left[\Theta_{{\zeta'}^{\beta}}'(g(w,z),\bar{g}(\zeta,\xi),
f(w,z))\right]_{
z=\xi+i\bar{\Theta}(\zeta,w,\xi)}\equiv_{w,\zeta,\xi}\\ &
\equiv \sum_{|\gamma| \leq |\beta|} \frac{
\underline{{\Cal L}}^{\gamma} \bar{f}(\zeta,\xi) \ \underline{u}_{\beta,\gamma}
((\underline{\Cal L}^{\delta} \bar{g}(\zeta,\xi))_{|\delta|\leq
|\beta|})}{
\underline{\Delta}(w,\zeta,\xi)^{2|\beta|-1}}\\
&
\equiv : \underline{\Omega}_{\beta}(w,\zeta,\xi,
(\partial_{\xi}^{\alpha^1}\partial_{\zeta}^{\gamma^1}\bar{h}(
\zeta,\xi))_{|\alpha^1|+ |\gamma^1| \leq
|\beta|})\\ & =: \underline{\omega}_\beta(w,\zeta,\xi),
\endaligned
\tag 5.3
$$
as formal power series in $(w,\zeta,\xi)$, where
$$
\aligned
&\underline{\Delta}(w,\zeta,\xi)=\underline{\Delta}(w,z,\zeta,\xi)|_{z=
\xi+i\bar{\Theta}(\zeta,w,\xi)}={\text{\rm det}}(\underline{\Cal L} \bar{g})=\\
&= {\text{\rm det}} \left(
\frac{\partial \bar{g}}{\partial \zeta}(\zeta,\xi)-i\Theta_{\zeta}(w,\zeta,z)
\frac{\partial \bar{g}}{\partial \xi}(\zeta,\xi)\right)|_{z=\xi+
i\bar{\Theta}(\zeta,w,\xi)}.
\endaligned
\tag 5.4
$$
\endproclaim

\remark{Remarks} 

1. As ${\text{\rm det}} (h_{i,\1_j})_{1\leq i,j\leq n}\neq 0$, then
${\text{\rm det}}(\underline{\Cal L}\bar{g}(0))\neq 0$ also, so
$\underline{\Delta}^{1-2|\beta|}\in \C\dl w,\zeta,\xi\dr$.

2. Putting $(\zeta, \xi)=(0,0)$, but after perharps shrinking $r$, we
first readily observe that $\underline{\Delta}^{1-2\vert \beta \vert }
(w,0,0) \in \Cal O ((r\Delta)^{m},\C)$ for some $r>0$, since
$\Theta_\zeta(w,0,0)
\in \Cal O ((r\Delta)^{m},\C)$ and since $\partial_\zeta^{\gamma^1}
\bar g(0,0)$ for $\vert \gamma^1\vert =1$ and $\partial_\xi^{\alpha^1} \bar g
(0,0)$ for $\vert \alpha^1\vert=1$ are constants, so that we observe
finally that
$$
\underline{\Omega}_\beta(w,0,0,(\partial_\xi^{\alpha^1} 
\partial_\zeta^{\gamma^1} \bar h
(0,0))_{
\vert \alpha^1 \vert +\vert \gamma^1 \vert \leq \vert \beta \vert}
)\in \Cal O((r\Delta)^{m},\C),
\tag 5.5
$$
for all $\beta \in \N_*^m$, because all the coefficients of the
$\underline{\Cal L}^\gamma$ belong to $\Cal
O((r\Delta)^{m},\C)$. Therefore, {\it the domains of convergence of
the $\underline{\omega}_\beta(w,0,0)$ are independent of $\beta$.}

3. $\underline{\Omega}_{\beta}$ arises after writing $\underline{\Cal
L}^{\gamma}\bar{h}(\zeta,\xi)$ as $\chi_{\delta}(w,z,\zeta,\xi,
(\partial_{\xi}^{\alpha^1}\partial_{\zeta}^{\gamma^1}\bar{h}(\zeta,\xi))_{|\alpha^1|+
|\gamma^1| \leq |\beta|})$ (by noticing that the coefficients of
$\underline{\Cal L}$ are analytic in $(w,z,\zeta,\xi)$) and by
replacing again $z$ by $\xi+i\bar{\Theta}(\zeta,w,\xi))$.
\endremark\medskip
 
\demo{Proof} 
Applying the $\underline{\Cal L}_j$, $1\leq j\leq m$, to eq. (5.1), we
obtain
$$
\underline{\Cal L}_j \bar{f} =-i\sum_{|\beta|=1} \underline{\Cal L}_j
\bar{g}^\beta \ \Theta_{{\zeta'}^{\beta}}'(g,\bar{g},f), \ \ \ \ \ 1\leq j\leq m,
\tag 5.6
$$
Then Cramer's rule applied to (5.6) yields immediately (5.3) for
$|\beta|=1$.

By induction, applying the $\underline{\Cal L}_j$, $1\leq j\leq m$, to
(5.3), we obtain
$$
\aligned
& {1\over \beta!}\sum_{|\beta_1|=1}
\underline{\Cal L}_j \bar{g}_{\lrcorner\beta_1} \Theta_{{\zeta'}^{\beta+\beta_1}}'
(g,\bar{g},f)\equiv\\ &
\equiv \sum_{|\gamma| \leq |\beta|} 
\frac{
\underline{\Cal L}^{\gamma+1_j} \bar{f} \ \underline{u}_{\beta,\gamma}
((\underline{\Cal L}^{\delta} \bar{g})_{|\delta|\leq |\beta|})}{
\underline{\Delta}^{2|\beta|-1}}+
\frac{
\underline{\Cal L}^{\gamma} \bar{f} \ \sum_{X_{\delta}} \frac{\partial 
\underline{u}_{\beta,\gamma}}{\partial X_{\delta}}
((\underline{\Cal L}^{\delta}\bar{g})_{|\delta|\leq |\beta|}) \
\underline{\Cal L}^{\delta+\1_j} \bar{g}}{
\underline{\Delta}^{2|\beta|-1}}+\\
& +\frac{
\underline{\Cal L}^{\gamma}\bar{f} \ \ \underline{u}_{\beta,\gamma}((
\underline{\Cal L}^{\delta}\bar{g})_{|\delta|\leq |\beta|}) (1-2|\beta|) 
\underline{\Cal L}_j(\underline{\Delta})}{
\underline{\Delta}^{2|\beta|}},
\endaligned
\tag 5.7
$$
where $\lrcorner\beta$ is the integer $k$ such that $\beta_k=1$ if
$|\beta|=1$ and where $((X_\delta)_{\vert \delta \vert \leq
\vert \beta \vert})$ denote indeterminates in place of
$((\underline{\Cal L}^{\delta}\bar{g})_{|\delta|\leq |\beta|})$ for
the functions $u_{\beta,\gamma}= u_{\beta,\gamma}((X_\delta)_{\vert
\delta \vert \leq \vert \beta \vert})$.

Now, applying Cramer's rule to (5.7), we get (5.3) with
$\beta+\beta_1$ instead of $\beta$, for all $|\beta_1|=1$, which
completes the proof of Proposition 5.2.
\qed\enddemo

\remark{Remark} 
Analogously, chosing instead the derivation $\Cal L^\beta$, we have
$$
\aligned
&
\theta_{\beta}'(\zeta,w,z):=\left[
\bar{\Theta}_{w^{\beta}}'(\bar{g}(\zeta,\xi),g(w,z), 
\bar{f}(\zeta,\xi))\right]_{\xi=
i\Theta(w,\zeta,z)} \equiv_{\zeta,w,z}
\equiv \\
&
\sum_{|\gamma|\leq |\beta|}
\frac{
{\Cal L}^{\gamma} f(w,z) u_{\beta,\gamma} (({\Cal L}^{\delta}
g(w,z))_{|\delta| \leq |\beta|})}{
\Delta(\zeta,w,z)^{2|\beta|-1}}\\
& =\Omega_{\beta}(\zeta,w,z,(\partial_z^{
\alpha}\partial_w^{\gamma}h(w,z))_{
|\alpha|+|\gamma|\leq |\beta|})\\
\endaligned
\tag 5.8
$$
as formal power series in $(\zeta, w,z)$, where
$$
\aligned
&
\Delta(\zeta,w,z)=\Delta(\zeta,\xi,w,z)|_{\xi=
z-i\Theta(w,\zeta,z)}={\text{\rm det}}({\Cal L}g)=\\ & {\text{\rm
det}} \left(
\frac{\partial g}{\partial w}(w,z)+i\bar{\Theta}_w(\zeta,w,\xi) 
\frac{\partial g}{\partial z}(w,z)\right)_{\xi=z-i\Theta(w,\zeta,z)}.
\endaligned
\tag 5.9
$$
\endremark\medskip
 
But these new identities offer no real new information, because

\proclaim{Lemma 5.10} 
We have

$\underline{\overline{\theta}}_{\beta}'(\zeta,w,z)\equiv
\theta_{\beta}'(\zeta,w,z)$,

$\underline{\overline{\Delta}}(\zeta,w,z)\equiv \Delta(\zeta,w,z)$,

$\overline{\underline{\Omega}}_{\beta}(\zeta,w,z,
(\partial_z^{\alpha}\partial_w^{\gamma}h(w,z))_{ |\alpha|+|\gamma|\leq
|\beta|})
\equiv
\Omega_{\beta}(\zeta,w,z,(\partial_z^{\alpha}\partial_w^{\gamma}h(w,z))_{
|\alpha|+|\gamma|\leq |\beta|})$,

$\underline{\overline{u}}_{\beta,\gamma}(({\Cal
L}^{\delta}g(w,z))_{|\delta|\leq |\beta|})\equiv u_{\beta,\gamma}
(({\Cal L}^{\delta} g(w,z))_{|\delta|\leq |\beta|})$.
\endproclaim

\demo{Proof} When $\zeta= \bar{w}$ and $\xi=\bar{z}$, the equations 
(5.8) and (5.9) are just conjugates of (5.3) and (5.4) respectively.
\qed\enddemo

We now give a characterization of S-nondegeneracy of a $\Cal C^\omega$
generic $M'$.

\proclaim{Proposition 5.11}
Assume that $M'$ is given in normal coordinates. Then $M'$ is
S-nondegenerate if and only if there exist
$\beta_1,\ldots,\beta_{m'}\in \N_*^{m'}$ and integers $l_1,\ldots,
l_{m'}$, $1\leq l_i \leq d'$, such that
$$
{\text{\rm det}} \left(\frac{\partial
{\Theta'}_{\zeta^{\beta_i}}^{l_i}}{\partial {w'}_j} (w', 0,
0)\right)_{1\leq i,j\leq m'}
\not\equiv_{w'} 0 \ \ \ \ \ \hbox{in} \ \C\dl w'\dr. 
\tag 5.12
$$
\endproclaim

\demo{Proof}
This condition can be seen after replacing $t'$ by $(w',
i\bar{\Theta}' (0, w', 0))= (w', 0)$ in eq. (1.1.6) of the
introduction and $\rho'(t', \tau')$ by $\xi'-z' + i\Theta'(w', \zeta',
z')$. Indeed, if $h=(g,f)=\hbox{Id}=(w',z')$, we have
$\rho'(h(w',0),\bar{h}(0))=f(w',0)$, so in the determinant (1.1.6), we
can already include the terms $\rho'(h(w',0),\bar{h}(0))=f(w',0)$
amongst all the equations $R_{\gamma}' $, and these terms satisfy the
nice minor property $\frac{\partial \rho'}{\partial z'}
=\hbox{Id}_{d'\times d'}$. Thus, it only remains to find equations
for $g(w',0)$ similar to (1.1.6). But the set
$\{R_{\beta'}'=0\}_{\beta' \in \N_*^{m'}}$ coincides with the set
$\{\Theta_{{\zeta'}^{\beta'}}'-\underline\Omega_{\beta'}=0\}_{\beta'
\in \N_*^{m'}}$ in eq. (5.3). Since $g(w',0)=w'$, then condition
(1.1.7) reduces exactly to (5.12).
\qed\enddemo

\proclaim{Proposition 5.13}
If $n=n'$, if $M'$ is S-nondegenerate at $p'$ and if $h$ is of formal
rank $n$, then the formal map $h\: (M, p)\to_{\Cal F} (M', p')$ is
S-nondegenerate.
\endproclaim

\demo{Proof}
This is immediate, after using the characterization given in
Proposition 5.11, the equations given in eq. (5.3) and the fact that
$\hbox{det}\left({\partial g_j\over \partial w_k}_{1\leq j,k\leq
m}(0,0)\right)\neq 0$.
\qed\enddemo

\head \S6. Proof of Proposition 1.5.2 (ii) \endhead

The purpose of this paragraph is to prove Proposition 1.5.2 (ii),
namely to establish

\proclaim{Proposition 6.1}
If $m\geq m'$, if $M'$ is S-nondegenerate at $p'$ and if $h$ induces a
formal map $(S_{\bar{p}}, p) \to_{\Cal F} (S_{\bar{p}'}, p')$ of
generic rank equal to $m'={\text{\rm dim}}_{\C} S_{\bar{p}'}$, then
the formal map $h\: (M, p)\to_{\Cal F} (M', p')$ is S-nondegenerate.
\endproclaim

\demo{Proof} Thus, assume that 
the formal map $h\: (S_{\bar{p}}, p) \to_{\Cal F} (S_{\bar{p}'}', p')$
is of generic rank equal to $m'$ (hence in particular, $m$ must
satisfy $m\geq m')$. After renumbering the {\it normal} coordinates
$(w',z')$ for $M'$ near 0, we claim that can assume that the minor
$\hbox{det} (\underline{\Cal L} \bar{g} (\zeta,
\xi))^{\square}: ={\text{\rm det}}\left(
\underline{\Cal L}_j \bar{g}_k(\zeta, \xi) \right)_{1\leq j,k\leq m'}$
satisfies ${\text{\rm det}}(\underline{\Cal L} \bar{g} (\zeta,
0))^{\square}\not\equiv 0$ in $\C\dl \zeta \dr$, {\it i.e.} more
precisely, that $[{\text{\rm det}}(\underline{\Cal L} \bar{g} (\zeta,
\xi))^{\square}]_{w=0,z=0,\xi=0}\not\equiv_\zeta 0$ in $\C\dl \zeta
\dr$. Indeed, for any power series $a(\zeta)$, then
$\underline{\Cal L}^{\gamma} a(\zeta)|_{w=z=\xi=0}=
(\partial^{|\gamma|} / \partial \zeta^{\gamma}) a(\zeta)$, because
$\underline{\Cal L}=\partial /\partial \zeta$ at $(w, z,\zeta,
\xi)= (0,0,\zeta,0)$. This property can be again readily seen by 
noticing that the induced mapping $h\vert_{S_0} \: (S_0, 0)\to_{\Cal
F} (S_0',0)$ is simply given by $(\zeta,0) \mapsto_{\Cal F}
(g(\zeta),0)$.

Now, it is easy to observe after using the adjoint matrix of
$(\underline{\Cal L} \bar{g} (\zeta,\xi))^{\square}$ (instead of its
inverse in case it is invertible) that the same calculation as the one
that we did in Proposition 5.2 yields immediately
$$
\aligned
& (\hbox{det} (\underline{\Cal L} \bar{g} (\zeta,
\xi))^{\square})^{2|\beta|-1}
\left[\Theta_{{\zeta'}^{\beta}}'(g(w,z),\bar{g}(\zeta,\xi), f(w,z))\right]_{
z=\xi+i\bar{\Theta}(\zeta,w,\xi)}\equiv_{w,\zeta,\xi}\\ &
\equiv \sum_{|\gamma| \leq |\beta|}
\underline{{\Cal L}}^{\gamma} \bar{f}(\zeta,\xi) \ \underline{u}_{\beta,\gamma}
((\underline{\Cal L}^{\delta} \bar{g}(\zeta,\xi))_{|\delta|\leq
|\beta|})=:
\underline{\Cal L}^{\gamma^{\sharp}(2|\beta|-1)}
\underline{\Omega}_{\beta}^1 / (\underline{\Cal L}^{\gamma^{\sharp}(2
|\beta|-1)} \Delta^{2\vert \beta \vert +1})
\endaligned
\tag 6.2
$$
Now the determinant $\underline{\Delta}(w,\zeta,\xi) :={\text{\rm
det}}(\underline{\Cal L}\bar{g} (\zeta, \xi))^{\square}$ satisfies
$\underline{\Delta}(0,\zeta,0)\not\equiv 0$ by assumption. This holds
if and only if there exists a multiindex $\gamma_{\sharp}\in \N_*^m$
such that $[\underline{\Cal L}^{\gamma_{\sharp}}
\underline{\Delta}(w,\zeta,\xi)]|_{w=z=\xi=0}
\neq 0$. Indeed, 
$\underline{\Cal L}^{\gamma} a(\zeta)|_{\zeta=0}= (\partial^{|\gamma|}
/ \partial \zeta^{\gamma}) a(\zeta)|_{\zeta=0}$ because
$\underline{\Cal L}=\partial /\partial \zeta$ at $(w, \zeta, z)=
(0,\zeta,0)$.

Moreover, we can choose $\gamma^{\sharp}$ to be minimal with respect
to the lexicographic order.

By applying now the derivation $(\underline{\Cal
L}^{\gamma_{\sharp}})^{2|\beta|-1}$ to eq. (6.2), a derivation which
satisfies $ (\underline{\Cal L}^{\gamma_{\sharp}})^{2|\beta|-1}
(\underline{\Delta}(w,\zeta,\xi)^{2|\beta|-1}) |_{w=z=\zeta=\xi=0}
\neq 0$, but 
$$
[\underline{\Cal L}^{\gamma}
(\underline{\Delta}(w,\zeta,\xi)^{2|\beta|-1})] |_{w=\zeta=\xi=0}=0, \
\
\forall \ \gamma < (\gamma_{\sharp})(2|\beta|-1)
\tag 6.3
$$
(since $\gamma^{\sharp}$ is minimal with respect to the lexicographic
order), we get that there exists a holomorphic remainder $C_{\beta}$
such that we can write
$$
\left[\Theta_{{\zeta'}^{\beta}}'(g(w,z),\bar{g}(\zeta,\xi), f(w,z))\right]_{
z=\xi+i\bar{\Theta}(\zeta,w,\xi)}+C_{\beta}(g(w,z), f(w,z),
\nabla^{|\gamma_{\sharp}|(2|\beta|-1)}\bar{g}
\tag 6.4
$$
$$\bar{g}(\zeta, \xi)) =
\underline{\Xi}_{\beta}(
w, \zeta, \xi, \nabla^{|\beta|+|\gamma_{\sharp}|(2|\beta|-1)}
\bar{h}(\zeta, \xi))=
\underline{\Cal L}^{\gamma^{\sharp}(2|\beta|-1)}
\underline{\Omega}_{\beta}^1 / (\underline{\Cal L}^{\gamma^{\sharp}(2
|\beta|-1)} \Delta)
$$
a remainder which has the property
$$
C_{\beta}(g(w,0), f(w,0),
\nabla^{|\gamma_{\sharp}|(2|\beta|-1)}\bar{g}(0,0))\equiv 0.
\tag 6.5
$$

In summary, amongst the equations ${R'}_{\beta}^{l}(w,0, 0, h(w,0),
\nabla^{|\beta|}\bar{h}(0))\equiv 0$ for $\beta\neq 0$, we have
obtained equations of the form
$$
\Theta_{{\zeta'}^{\beta}}'(g(w,0),\bar{g}(0), 0)=
\underline{\Xi}_{\beta}(
w,0,0, \nabla^{|\beta|+|\gamma_{\sharp}|(2|\beta|-1)}
\bar{h}(0,0))
\tag 6.6
$$
(recall $f(w,0)\equiv 0$). Then condition (1.1.6) defining
S-nondegeneracy is clearly satisfied, because (5.12) holds and
$\hbox{det} ({\Cal L} g(w,0))^{\square}\not\equiv 0$.
\qed\enddemo

\heading \S7. Propagation of analyticity along Segre chains \endheading

Let us recall that we started in the introdution with the family of
equations
$$
0=\underline{\Cal L}^{\gamma} \rho'(h(t), \bar{h}(\tau)) :=
R_{\gamma}'(t,\tau, h(t), \nabla^{|\gamma|} \bar{h}(\tau))=0
\tag 7.1
$$
satisfied by $h(t)$ as $\rho(t,\tau)=0$.

To treat S-nondegenerate maps, we delineate the following statement.

\proclaim{Theorem 7.2}
Let $h\: (M, p)\to_{\Cal F} (M', p')$ be a formal holomorphic map with
$M$ given by equations $(2.1.1)$ and minimal at $p$. Assume that there
exist $\kappa_0 \in \N_*$ and a finite family of relations of the form
$$
{\Cal X}_{\lambda}' (t,\tau, h(t), \nabla^{\kappa_0}
\bar{h}(\tau))\equiv 0 \ \ \ \ \ {\text{\rm on}} \ \ {\Cal
M}=\{\rho(t,\tau)=0\},
\tag 7.3
$$
for $\lambda= 1, \ldots, \Lambda$, $\Lambda \in \N_*$, the ${\Cal
X}_{\lambda}'= {\Cal X}_{\lambda}'(t, \tau, t', \nabla^{\kappa_0})$
being holomorphic from a neighborhood of $0\times 0\times 0\times
\nabla^{\kappa_0} \bar{h}(0)$ in 
$\C^n\times \C^n\times \C^{n'} \times \C^{n'N_{n,\kappa_0}}$, and
$\C$-valued, such that the following nondegeneracy condition holds:

$(*)$ There exist $\lambda_1, \ldots, \lambda_{n'} \in \dl 1, \Lambda
\dr$ such that
$$
{\text{\rm det}} \left(
\frac{\partial {\Cal X}_{\lambda_j}'}{\partial t_k'}(w, 
i\bar{\Theta}(0, w, 0), 0, h(w, i\bar{\Theta}(0, w,0)),
\nabla^{\kappa_0}\bar{h}(0))
\right)_{1\leq j, k\leq n'} \not\equiv_w 0 \ \ \hbox{in} \ \ \C \dl w \dr.
\tag 7.4
$$
Then $h$ is convergent.
\endproclaim

To treat S-finite maps and in the same time S-solvable maps, we
delineate the following statement.

\proclaim{Theorem 7.5}
{\text{\rm [BER99]}} Let $h\: (M, p)\to_{\Cal F} (M', p')$ be a formal
holomorphic map, with $M$ given by equation (2.1.1) and minimal at
$p$. Assume that

$(**)$ There exist $\kappa_0\in \N_*$, $N_j'\in \N_*$, $1\leq j\leq
n'$ and monic Weierstrass polynomials of the form $P_j'(t, \tau, t_j',
\nabla^{\kappa_0})$,
$$
P_j'(t, \tau, t_j', \nabla^{\kappa_0}) = {t'}_j^{N_j'} + \sum_{1\leq
k\leq N_j'} A_{j,k}'(t, \tau,
\nabla^{\kappa_0}){t'}_j^{N_j'-k},
\tag 7.6
$$
the $A_{j,k}'$ being holomorphic from a neighborhood of $0\times 0
\times \nabla^{\kappa_0} \bar{h}(0)$ in $\C^n\times \C^n \times\C^{n'
N_{n,\kappa_0}}$, $1\leq j\leq n'$, and $\C$-valued, such that the
following formal identities
$$
P_j'(t, \tau, h_j(t), \nabla^{\kappa_0} \bar{h}(\tau))\equiv 0, \ \ \
\ \ 1\leq j\leq n',
\tag 7.7
$$
hold if $\rho(t,\tau)=0$. Then $h$ is convergent.
\endproclaim

To apply Theorem 7.5 to S-finite formal maps, we shall consider a
transformation of the complex analytic subset of a neighborhood of
$0\times 0\times 0 \times \nabla^{\kappa_0} \bar{h}(0)$ in $\C^n\times
\C^n \times \C^{n'} \times \C^{n'N_{n,\kappa_0}}$ defined by the
equations
$$
\SS^1:= \{(t, \tau, t', \nabla^{\kappa_0}) \: 
\rho(t, \tau)=0, R_{\gamma}'(t, \tau, t', \nabla^{\kappa_0} \bar{h}(\tau))=0 \ 
\forall |\gamma| \leq \kappa_0 \},
\tag 7.8
$$
where $\kappa_0$ is large enough.

Thanks to the assumption that ${\text{\rm dim}}_{p'}(\V_p')=0$
($\V_p'= \SS^1 \cap (0\times 0\times \C^{n'} \times \nabla^{ \kappa_0}
\bar{h}(0))$), we can utilize a classical transformation of the
equations of $\SS^1$ which consists in replacing the $R_{\gamma}'$'s
by the canonical defining functions of the ramified analytic cover
$\pi\: S^1 \to \C^n\times \C^n \times 0\times \C^{n'N_{n,\kappa_0}}$
(Whitney's contruction, [Ch1], Chapter 1, paragraph 4, pp. 42-51) to
obtain:

\proclaim{Lemma 7.9}
If $h$ is S-finite at $p$, then there exist $N\in \N_*$, $\kappa_0\in
\N_*$ and Weierstrass polynomials $P_j'(t, \tau, t_j',
\nabla^{\kappa_0})= {t_j'}^N + \sum_{1\leq k\leq N} A_{j,k}(t,\tau,
\nabla^{\kappa_0}){t_j'}^{N-k}$, 
with $A_{j,k}$ being holomorphic in a neighborhood of $0\times 0\times
\nabla^{\kappa_0}\bar{h}(0)$ in $\C^n\times \C^n \times
\C^{n'N_{n,\kappa_0}}$, $1\leq j\leq n'$, such that 
$\SS^1$ is contained in the complex analytic set
$$
\E^1:=\{(t,\tau, t', \nabla^{\kappa_0})\:
P_j'(t,\tau, t', \nabla^{\kappa_0})=0\}
\tag 7.10
$$
and the following formal equalities
$$
P_j'(t,\tau, h_j(t), \nabla^{\kappa_0} \bar{h}(\tau))\equiv 0, \ \ \ \
\ 1\leq j\leq n',
\tag 7.11
$$
hold if $(t,\tau)$ satisfy $\rho(t,\tau)=0$.
\endproclaim

\demo{Proof}
Existence of $\E^1 \supset \SS^1$ follows by taking a suitable subset
of the set of canonical defining functions of the analytic cover
$\pi\: \SS^1 \to \C^n\times \C^n \times 0 \times
\C^{n'N_{n,\kappa_0}}$, namely the set of functions
which are polynomial in a single variable $t_j'$ (see [Ch1], p.42).

Let ${\Cal I}_{\SS^1}$ denote the ideal $(R_\gamma')_{|\gamma|\leq
\kappa_0}$ in ${\Cal O}_t\times {\Cal O}_{\tau} \times {\Cal O}_{t'}
\times {\Cal O}_{\nabla^{\kappa_0}}$. Each $P_j'$ vanishes over
$\SS^1$. Thus, by the Nullstellensatz, there exist integers $M_j\in
\N_*$ such that ${P'}_j^{M_j} \in {\Cal I}_{\SS^1}$, $1\leq j\leq
n'$. We deduce $P_j'(t,\tau, h(t), \nabla^{\kappa_0} \bar{h}
(\tau))^{M_j}\equiv 0$ when $\rho(t,\tau)=0$, whence $P_j'(t,\tau,
h(t), \nabla^{\kappa_0}\bar{h}(\tau))\equiv 0$ as desired.
\qed\enddemo

The main feature of S-finite maps is that after applying Lemma 7.9, we
see that each component $h_j(t)$ is entire over $\bar{h}(\tau)$ and
its jets, that is to say, it satisfies a monic polynomial equation
(7.11).

In particular, as S-solvable CR maps are S-finite, they satisfy Lemma
7.9 above. In truth, they satisfy a much more aesthetic relation.

\proclaim{Lemma 7.12}
If $h$ is S-solvable at $p\in M$, then there exist $N\in \N_*$,
$\kappa_0\in \N_*$, and analytic functions $A_j'(t,\tau,
\nabla^{\kappa_0})$ holomorphic in a neighborhood of $0\times 0\times
\nabla^{\kappa_0}\bar{h}(0)$ in $\C^n\times \C^n\times
\C^{n'N_{n,\kappa_0}}$, such that the following formal equalities
$$
h_j(t)\equiv A_j'(t,\tau, \nabla^{\kappa_0}\bar{h}(\tau)), \ \ \ \ \
1\leq j\leq n',
\tag 7.13
$$
hold if $(t,\tau)$ satisfy $\rho(t,\tau)=0$.
\endproclaim

\demo{Proof}
It suffices to apply the analytic implicit function theorem to the
collection of equations
$$
0= \underline{\Cal L}^{\gamma} \rho'(t', \bar{h}(\tau))=
R_{\gamma}'(t,\tau, t', \nabla^{|\gamma|}\bar{h}(\tau)), \ \ \ \ \
|\gamma|\leq \kappa_0,
\tag 7.14
$$
to solve them in $t'$ as $t_j'= A_j'(t,\tau,\nabla^{\kappa_0})$. In
fact, $h$ is S-solvable if and only if eq. (7.14) is solvable
(implicit function theorem) in $t'$ for $\kappa_0$ large enough.
\qed\enddemo

Heuristically we would like to say that this relation (7.13), compared
with the relation (7.7), shows {\it a posteriori} that in a certain
sense, a S-finite maps is {\it almost solvable} in terms $\bar h$ and
its jets $\nabla^\kappa \bar h$. {\it However, in the S-nondegenerate
case, the map $h$ is neither solvable nor entire over $\bar h$ and its
jets $\nabla^\kappa \bar h$, nor almost solvable in any sense.}

\head \S8. Proof of Theorem 7.2: Step I \endhead

This paragraph is devoted to the first step in the proof of Theorem
7.2.

\demo{Proof of Theorem 7.2}
First of all, we extract from the system of $({\Cal
X}_{\lambda}')_{1\leq \lambda \leq \Lambda}$ the subsystem $({\Cal
X}_{\lambda_j}')_{1\leq j\leq n'}$, we denote it by $({\Cal
X}_j')_{1\leq j\leq n'}$ and we show that if $h$ satisfies
$$ 
{\Cal X}_j'(t,\tau, h(t), \nabla^{\kappa_0} \bar{h}(\tau))\equiv 0 \ \
\ \ \
\hbox{on} \ \ {\Cal M}=\{\rho(t,\tau)=0\},
\tag 8.1
$$
for $j=1,\ldots, n'$, with ${\Cal X}_1',\ldots, {\Cal X}_{n'}'$
satisfying eq. (7.4), then $h$ is convergent.

We conduct the proof by induction on the length $k\in \N$ of a Segre
chain $\Gamma_{\underline{\Cal L} {\Cal L}}^k$ up to $k= 2\mu_p$,
where $\mu_p$ denotes the Segre type of ${\Cal M}$ at $p$ (see Theorem
3.2.2) and the induction processus $({\Cal I}_k) \Rightarrow ({\Cal
I}_{k+1})$ will be divided in two essential steps.

Here is the $k$-th induction assumption:

\smallskip

$({\Cal I}_k)$ There exists a $\C^{2n'}$-valued holomorphic map
$\Psi^k\: (\delta \Delta^m)^k \to \Delta^{2n'}$, $\delta >0$, such
that the composition
$$
h^c \circ \Gamma_{\underline{\Cal L}{\Cal L}}^k(w_1,\ldots, w_k)\equiv
\Psi^k (w_1,\ldots, w_k)\in \C\{w_1,\ldots,w_k\}^{2n'}
\tag 8.2
$$
is convergent for $w_1,\ldots,w_k\in \delta\Delta^m$.

\smallskip

It will be clear soon that $({\Cal I}_{2\mu_p})$ implies that $h^c$ is
convergent in a neighborhood of $p$ in $\Cal M$. Recall that
$h^c=(h,\bar{h})$, so eq. (8.2) above is a statement about the formal
complexified map $h^c \: ({\Cal M}, p^c) \to_{\Cal F} ({\Cal M}',
{p'}^c)$.

As announced above, the proof that $({\Cal I}_{k}) \Rightarrow ({\Cal
I}_{k+1})$ involves two essential steps:

\smallskip
\noindent
{\bf First step:} $[({\Cal I}_k)$ and $(*_{k-1})]
\Rightarrow (*_k)$, and 

\smallskip
\noindent
{\bf second step:} $(*_k)\Rightarrow ({\Cal I}_{k+1})$, where

\smallskip
$(*_k)$ For all $\kappa\in \N$ and all $\beta \in \N^n$ with $|\beta|
\leq \kappa$, there exist $\C^{2n'}$-valued holomorphic maps
$\Psi_{\beta}^k\: (\delta \Delta^m)^k \to \Delta^{2n'}$, $\delta >0$,
such that the composition
$$
(\nabla^{\kappa} h^c) \circ \Gamma_{\underline{\Cal L}{\Cal L}}^k
(w_1,\ldots, w_k)\equiv (\Psi_{\beta}^k (w_1,\ldots,
w_k))_{|\beta|\leq \kappa}
\in \C\{w_1,\ldots,w_k\}^{2n'N_{n,\kappa}}
\tag 8.3
$$
is convergent for $w_1,\ldots,w_k\in \delta\Delta^m$.

\remark{Remark} The positive
number $\delta$ may shrink to $0$ as $\kappa$ goes to $\infty$.
\endremark\medskip

The proof of the first step occupies this \S8 and the proof of step
two is postponed to \S9. Now, let us explain how we end-up the proof.

\proclaim{Assertion 8.4}
If $({\Cal I}_{2\mu_p})$ is satisfied, then $h\: (M, p) \to_{\Cal F}
(M', p')$ is convergent.
\endproclaim

\demo{Proof}
We know by Corollary 3.2.6 that for each $\delta >0$, there exist
$w_{(2\mu_p)}^* \in (\delta \Delta^m)^{2\mu_p}$ such that
$\Gamma_{\underline{\Cal L}{\Cal L}}^{2\mu_p} (w_{(2\mu_p)}^*) = p^c$
and such that the rank of $\Gamma_{\underline{\Cal L}{\Cal
L}}^{2\mu_p}$ at $w_{(2\mu_p)}^*$ equals ${\text{\rm dim}}_{\C} {\Cal
M}$. (Recall ${\Cal M}$ is $\{{\Cal L}, \underline{\Cal L}\}$-minimal
at $p^c$.) Thus $\Gamma_{\underline{\Cal L}{\Cal L}}^{2\mu_p}$ can
provide holomorphic coordinates over a small neighborhood of $p^c$ in
${\Cal M}$. By $({\Cal I}_{2\mu_p})$, we have that $h^c \circ
\Gamma_{\underline{\Cal L}{\Cal L}}^{2\mu_p}$ is convergent. This
proves that $h^c \: ({\Cal M}, p) \to_{\Cal F} ({\Cal M}', p')$ is
convergent. Thus $h$ is convergent.
\qed\enddemo

\remark{Remark}
According to Theorem 2.2.6, we know that $\hbox{gen-rk}_{\C}
(\Gamma_{\underline{\Cal L}{\Cal L}}^{\mu_p})=2m+d= {\text{\rm
dim}}_{\C} {\Cal M}$ already and we shall prove in particular that
$h\circ \Gamma_{\underline{\Cal L}{\Cal L}}^{\mu_p}$ is convergent
(this is $({\Cal I}_{\mu_p})$). Does it imply already that $h$
converges? Yes, and the result is due to Eakin and Harris [EH]: {\it
Let $a(y)$ be a formal power series and let $\varphi\: (\C_x^{\nu},0)
\to (\C_y^{\mu}, 0)$ be a local holomorphic map of generic rank
$\nu$. If $a\circ \varphi(x)\in \C\{x\}$ is convergent, then $a(y)$ is
convergent}. But the iteration process up to the $2\mu_p$-th step is
free here and we can conclude the end of the proof in a much more
elementary way after using Assertion 8.4 instead of applying the
theorem of Eakin and Harris$\ldots$!
\endremark\medskip

First of all, we prove $({\Cal I}_1)$ by applying Theorem 1.3.2 to the
family of equations
$$
{\Cal X}_j'(\Gamma_{\underline{\Cal L}{\Cal L}}^1(w_1), h\circ
\Gamma_{\underline{\Cal L}{\Cal L}}^1(w_1),
(\nabla^{\kappa_0} \bar{h}) \circ
\Gamma_{\underline{\Cal L}{\Cal L}}^1(w_1))\equiv 0, 
\ \ \ \ \ 1\leq j\leq n',
\tag 8.5
$$
or equivalently, since $\Gamma_{\underline{\Cal L}{\Cal L}}^1(w_1)=
{\Cal L}_{w_1}(0)= (w_1, i\bar{\Theta}(0, w_1, 0), 0, 0)$,
$$
{\Cal X}_j'(w_1, i\bar{\Theta}(0, w_1, 0), h(w_1, i\bar{\Theta}(0,
w_1, 0)),
\nabla^{\kappa_0}\bar{h}(0))\equiv 0, \ \ \ \ \ 
1\leq j\leq n'.
\tag 8.6
$$
The main assumption (7.4) of Theorem 7.2 fits the hypothesis of
Theorem 1.3.2 exactly. We thus get that the series $h(w_1,
i\bar{\Theta}(0, w_1, 0))= h\circ \Gamma_{\underline{\Cal L}{\Cal
L}}^1(w_1)\equiv \Psi^1(w_1)\in \C\{w_1\}$ is convergent. On the other
hand, $\bar{h} \circ \Gamma_{\underline{\Cal L}{\Cal L}}^1(w_1)
=\bar{h}(0)$ is clearly convergent, which completes our checking of
$({\Cal I}_1)$.

\smallskip

To prove that $[({\Cal I}_k) \ \hbox{and} \ (*)_{k-1}] \Rightarrow
(*_k)$, we first make preliminary considerations.

\smallskip

Let us introduce two systems of $d$-vector fields $\Upsilon$ and
$\underline{\Upsilon}$ which are both complementary to the subsystem
$\{{\Cal L}, \underline{\Cal L}\}$ (in order that each system $\{{\Cal
L},
\underline{\Cal L}, \Upsilon\}$ and $\{{\Cal L}, \underline{\Cal L},
\underline{\Upsilon}\}$ spans $T{\Cal M}$) and which can be written
$$
\Upsilon= \frac{\partial }{\partial z}
+(1-i\Theta_z(w, \zeta, z)) \frac{\partial }{\partial \xi} \ \ \ \ \
\hbox{and} \ \ \ \ \
\underline{\Upsilon}= \frac{\partial }{\partial \xi}+
(1+i\bar{\Theta}_{\xi} (\zeta, w, \xi)) \frac{\partial }{\partial z}.
\tag 8.7
$$
We observe immediately that ${\Cal L}$ and $\underline{\Upsilon}$
commute and that $\underline{\Cal L}$ and $\Upsilon$ commute also. Of
course, this commutation property can be written in terms of their
flows:
$$
{\Cal L}_{w_1} \circ \underline{\Upsilon}_{\xi_1}(p)=
\underline{\Upsilon}_{\xi_1} \circ {\Cal L}_{w_1}(p) \ \ \ \ 
\hbox{and} \ \ \ \ \underline{\Cal L}_{\zeta_1} \circ 
\Upsilon_{z_1} (p) = \Upsilon_{z_1} \circ \underline{\Cal L}_{\zeta_1}(p),
\tag 8.8
$$
for all $w_1, \zeta_1 \in \delta\Delta^m$, $z_1, \xi_1\in
\delta\Delta^d$, $p\in {\Cal M}$ and $\delta >0$ small. Let $q(x)$
denote a formal $2n$-vectorial power series in $\C\dl x\dr^{2n}$
vanishing at $0$, where $x\in \C^{\nu}$ and $\nu\in \N_*$.

From this observation and from Lemma 4.5, we deduce,
$$
\aligned
& h\circ \Upsilon_{z_1} \circ \underline{\Cal L}_{\zeta_1} (q(x))
\equiv h\circ \underline{\Cal L}_{\zeta_1} \circ \Upsilon_{z_1} (q(x))
\equiv h\circ \Upsilon_{z_1}(q(x)),\\ &
\bar{h}\circ \underline{\Upsilon}_{\xi_1} \circ {\Cal L}_{w_1} (q(x)) \equiv 
\bar{h} \circ {\Cal L}_{w_1} \circ \underline{\Upsilon}_{\xi_1} (q(x)) \equiv
\bar{h} \circ \underline{\Upsilon}_{\xi_1} (q(x)).
\endaligned
\tag 8.9
$$
More generally and similarly
$$
\aligned
& (\nabla^{\kappa} h)\circ \Upsilon_{z_1} \circ \underline{\Cal
L}_{\zeta_1} (q(x)) \equiv (\nabla^{\kappa} h) \circ
\Upsilon_{z_1} (q(x)),\\
& (\nabla^{\kappa} \bar{h})
\circ \underline{\Upsilon}_{\xi_1} \circ {\Cal L}_{w_1} (q(x)) \equiv
(\nabla^{\kappa}\bar{h}) \circ \underline{\Upsilon}_{\xi_1}(q(x)).
\endaligned
\tag 8.10
$$

Now, we remark that the fundamental identity (7.4) also comes together
with a conjugate identity
$$
\bar{\Cal X}_j'(\tau, t, \bar{h}(\tau), \nabla^{\kappa_0} h(t)) \equiv 0 
\ \ \ \ \ \hbox{on} \ \ {\Cal M} = \{\rho(t,\tau)=0\},
\tag 8.11
$$
for $j= 1,\ldots, n'$.

By the way, this second identity can be simply derived from (7.4) by
specifying $\tau= \bar{t}$ in (7.4), by conjugating (7.4), which
yields $\bar{\Cal X}_j'(\bar{t}, t, \bar{h}(\bar{t}),
\nabla^{\kappa_0} h(t))\equiv 0$ on $M=\{\rho(t, \bar{t})=0\}$ and
then complexifying $(\bar{t})^c=\tau$, using that $M\subset {\Cal M}$
is a maximally real manifold. By this, we achieve a duplication, in
accordance with the heuristic principle that there is no privilegiate
choice between $h$ and $\bar{h}$.

This second identity will be used to achieve the second inductive step
$(*_k)\Rightarrow ({\Cal I}_{k+1})$ in case $k$ is odd and to achieve
the first inductive step $[({\Cal I}_k)$ and $(*_{k-1})] \Rightarrow
(*_k)$ in case $k$ is even.

We shall only prove $[({\Cal I}_{k})$ and $(*_{k-1})] \Rightarrow
(*_k)$ and $(*_k)\Rightarrow ({\Cal I}_{k+1})$ in case $k$ is odd. The
even case is completely similar. We would like to mention that the
scheme of our proof differs in an essential way from the proofs given
in [BER97] [BER99]. Whereas a solution $h_j(t)=A_j'(t,\tau,
\nabla^{\kappa_0}\bar{h}(\tau))$ or $P_j'(t, \tau, h_j(t),
\nabla^{\kappa_0} \bar{h}(\tau))=0$ (see \S11) provides an immediate explicit
iteration processus along Segre chains, we cannot here iterate
directly $h$ in terms of $\nabla^{\kappa_0}\bar{h}$ along Segre
chains. But we can use step by step the powerfulness of Theorem 1.3.2
and we shall have to check step by step the convergence of jets of $h$
without such an explicit relation of solvability.

As $k$ is odd, we write $k=2l+1$, $l\in \N$, and develope the
expression of $\Gamma_{\underline{\Cal L}{\Cal L}}^{2l+1}$ in the long
explicit form (just by definition)
$$
\Gamma_{\underline{\Cal L}{\Cal L}}^{2l+1} (w_{l+1}, \zeta_l, w_l, \ldots, 
\zeta_1, w_1) = 
{\Cal L}_{w_{l+1}} \circ \underline{\Cal L}_{\zeta_l} \circ {\Cal
L}_{w_l} \circ \cdots \circ
\underline{\Cal L}_{\zeta_1} \circ {\Cal L}_{w_1} (0).
\tag 8.12
$$
This expression in terms of vector fields is geometric, invariant and
appropriate for our understanding of the sequel, especially for the
calculation of high order derivatives. We shall therefore keep such
long explicit forms in the formalism to perform the calculations
below. Thus, the remainder of \S8 is devoted to the proof of $[{\Cal
I}_{2l+1} \ \hbox{and} \ (*_{2l}) ] \Rightarrow (*_{2l+1})$.

By substituting the series
$$
\underline{\Upsilon}_{\xi_1} \circ 
\Gamma_{\underline{\Cal L}{\Cal L}}^{2l+1} (w_{l+1}, \zeta_l, w_l, \ldots, 
\zeta_1, w_1)\in \C\dl w_1, \zeta_1, \ldots, 
w_l, \zeta_l, w_{l+1}\dr^{2n}
\tag 8.13
$$ for the variables $(t,\tau)\in {\Cal M}$ in the fundamental
identities (7.3), we can read (7.3) as
$$
{\Cal X}_j'(\underline{\Upsilon}_{\xi_1} \circ {\Cal L}_{w_{l+1}}
\circ \cdots \circ
\underline{\Cal L}_{\zeta_1} \circ {\Cal L}_{w_1} (0), \
h \circ \underline{\Upsilon}_{\xi_1} \circ {\Cal L}_{w_{l+1}}
\circ \underline{\Cal L}_{\zeta_l} \circ {\Cal L}_{w_l} \circ \cdots \circ
\underline{\Cal L}_{\zeta_1} \circ {\Cal L}_{w_1} (0), 
$$
$$ \ \ \ \ \ \ \ \ \ \ \ \ \ \ \ \ (\nabla^{\kappa_0} \bar{h}) \circ
\underline{\Upsilon}_{\xi_1} \circ {\Cal L}_{w_{l+1}}
\circ \underline{\Cal L}_{\zeta_l} \circ {\Cal L}_{w_l} \circ \cdots \circ
\underline{\Cal L}_{\zeta_1} \circ {\Cal L}_{w_1} (0))\equiv 0,
\ \ \ \ \ 1\leq j\leq n'
\tag 8.14
$$
or after a crucial simplification of the last term:
$$
{\Cal X}_j'(\underline{\Upsilon}_{\xi_1} \circ {\Cal L}_{w_{l+1}}
\circ \cdots \circ
\underline{\Cal L}_{\zeta_1} \circ {\Cal L}_{w_1} (0), \
h \circ \underline{\Upsilon}_{\xi_1} \circ {\Cal L}_{w_{l+1}}
\circ \underline{\Cal L}_{\zeta_l} \circ {\Cal L}_{w_l} \circ \cdots \circ
\underline{\Cal L}_{\zeta_1} \circ {\Cal L}_{w_1} (0), 
$$
$$ \ \ \ \ \ \ \ \ \ \ \ \ \ \ \ \ (\nabla^{\kappa_0} \bar{h}) \circ
\underline{\Upsilon}_{\xi_1}
\circ \underline{\Cal L}_{\zeta_l} \circ {\Cal L}_{w_l} \circ \cdots \circ
\underline{\Cal L}_{\zeta_1} \circ {\Cal L}_{w_1} (0))\equiv 0,
\ \ \ \ \ 1\leq j\leq n'
\tag 8.15
$$
in view of identity (8.10). The goal is to prove that the jets of $h$
on the $2l+1$-th Segre chain all converge, {\it i.e.} that for all
$\kappa\in\N$, then $(\nabla^{\kappa} h) \circ \Gamma_{\underline{\Cal
L} {\Cal L}}^{2l+1} (w_{l+1}, \zeta_l, w_l, \ldots,
\zeta_1, w_1)$
is convergent. To this aim, we claim that it suffices to prove that
$$
({\Cal L}^{\gamma} \underline{\Upsilon}^{\delta}h ) \circ
\Gamma_{\underline{\Cal L}{\Cal L}}^{2l+1}(w_{l+1}, \zeta_l, w_l, \ldots, 
\zeta_1, w_1) \in \C\{w_1, \zeta_1, \ldots, w_l, \zeta_l, w_{l+1}\}
\tag 8.16
$$
for all $\gamma \in \N^m$, $\delta \in \N^d$, with $|\gamma|+|\delta|
\leq \kappa$. In fact, it is easy to see that 
property (8.16) is equivalent to all $(\nabla^{\kappa} h) \circ
\Gamma_{\underline{\Cal L} {\Cal L}}^{2l+1} (w_1, \zeta_1, \ldots,
w_l, \zeta_l, w_{l+1})$ being convergent, since $\{{\Cal L},
\underline{\Upsilon}\}$ spans the $(w,z)$-space and has analytic
coefficients.

Recall that by definition of the flow, for a power series, $a\: {\Cal
M} \to \C$, $a(p)=0$, one has $\frac{\partial }{\partial w_1}(a\circ
{\Cal L}_{w_1} (q(x)))= ({\Cal L} a) \circ {\Cal L}_{w_1} (q(x))$ and
$\frac{\partial }{\partial \xi_1}|_{\xi_1=0} (a\circ
\underline{\Upsilon}_{\xi_1} (q(x)))=
(\underline{\Upsilon} a)(q(x))$ (in symbolic notations, omitting
indices for $w_{1,1},\ldots, w_{1,m}$ and
$\xi_{1,1},\ldots,\xi_{1,d}$).

Of course, (8.16) is satisfied for $\gamma=0$ and $\delta=0$, since we
have assumed that $(\Cal I_{2l+1})$ holds true.

Now, assume by induction that (8.16) is satisfied for all $\gamma,
\delta$ with $|\gamma|+|\delta|\leq \lambda$, some $\lambda \in \N_*$.
Applying to equation (8.14) the derivatives
$$
\left[\frac{\partial^{|\gamma|+1} }{\partial w_{l+1}^{\gamma+\1_r^m}}
\frac{\partial^{|\delta|} }{\partial
\xi_1^{\delta}} \ \ \bullet \ \ \right]_{\xi_1=0}\ \ \ \ \ \hbox{and} \ \ \
\ \
\left[\frac{\partial^{|\gamma|} }{\partial w_{l+1}^{\gamma}}
\frac{\partial^{|\delta|+1} }{\partial \xi_1^{\delta+\1_s^d}}
\ \ \bullet \ \ 
\right]_{\xi_1=0}
\tag 8.17
$$
where $\1_r^m=(0,\ldots, 1,\ldots, 0)\in \N^m$ with $1$ at the $r$-th
place, and where $\1_s^d= (0,\ldots, 1,\ldots, 0)$, with $1$ at the
$s$-th place, we will obtain after applying the induction assumption
that there exist two families of analytic series $A_{j, \gamma,
\delta,
\1_r^m}$ and $B_{j, \gamma, \delta, \1_s^d}$ such that
$$
\sum_{k=1}^{n'} \frac{\partial {\Cal X}_j'}{\partial t_k'}
( {\Cal L}_{w_{l+1}} \circ \cdots \circ {\Cal L}_{w_1} (0), h \circ
{\Cal L}_{w_{l+1}} \circ \cdots \circ {\Cal L}_{w_1}(0),
\tag 8.18
$$
$$
(\nabla^{\kappa_0} \bar{h}) \circ
\underline{\Cal L}_{\zeta_l} \circ {\Cal L}_{w_l} \circ \cdots \circ
{\Cal L}_{w_1}(0))({\Cal L}^{\gamma+ \1_r^m}
\underline{\Upsilon}^{\delta} h_k)
\circ {\Cal L}_{w_{l+1}} \circ \cdots \circ
{\Cal L}_{w_1} (0)+
$$
$$ \ \ \ \ \ \ \ \ \ \ \ \ + A_{j,\gamma, \delta, \1_r^m} (w_1,
\zeta_1, \ldots, w_l, \zeta_l, w_{l+1})\equiv 0, \ \ \ \ 1\leq j\leq
n'
$$
and
$$
\sum_{k=1}^{n'} \frac{\partial {\Cal X}_j'}{\partial t_k'}
( {\Cal L}_{w_{l+1}} \circ \cdots \circ {\Cal L}_{w_1} (0), h \circ
{\Cal L}_{w_{l+1}} \circ \cdots \circ {\Cal L}_{w_1}(0),
\tag 8.19
$$
$$
(\nabla^{\kappa_0} \bar{h}) \circ
\underline{\Cal L}_{\zeta_l} \circ {\Cal L}_{w_l} \circ \cdots \circ
{\Cal L}_{w_1}(0))({\Cal L}^{\gamma}
\underline{\Upsilon}^{\delta+\1_s^d} h_k)
\circ {\Cal L}_{w_{l+1}} \circ \cdots \circ
{\Cal L}_{w_1} (0) +$$
$$ \ \ \ \ \ \ \ \ \ \ \ \ + B_{j,\gamma, \delta, \1_s^d} (w_1,
\zeta_1, \ldots, w_l, \zeta_l, w_{l+1})\equiv 0, \ \ \ \ 1\leq j\leq
n'.
$$
Indeed, all the terms appearing in $A_{j, \gamma, \delta, \1_r^m}$ (or
$B_{j,\gamma, \delta, \1_s^d}$) involve three sorts of terms.

Firstly, they involve derivatives of $h$ of the form $({\Cal
L}^{\alpha} \underline{\Upsilon}^{\beta} h)\circ {\Cal L}_{w_{l+1}}
\circ \cdots \circ {\Cal L}_{w_1}(0)$ with $|\alpha|+|\beta| \leq
|\gamma| +|\delta|$, which are already analytic, by our induction
assumption on $\lambda=|\gamma| +|\delta|$, together with derivatives
of the form
$$
\sum_{|\alpha_1|+|\beta_1|+|\gamma_1|+|\delta_1| \leq
|\gamma|+|\delta|+1}
\frac{\partial^{|\alpha_1|+|\beta_1|+|\gamma_1|+|\delta_1|} {\Cal X}_j' }{
\partial t^{\alpha_1} \partial\tau^{\beta_1} \partial {t'}^{\gamma_1} 
(\partial \nabla^{\kappa_0})^{\delta_1}} ({\Cal L}_{w_{l+1}} \circ
\cdots \circ {\Cal L}_{w_1}(0),
\tag 8.20
$$
$$
h\circ {\Cal L}_{w_{l+1}} \circ \cdots \circ {\Cal L}_{w_1}(0), \
(\nabla^{\kappa_0} \bar{h}) \circ \underline{\Cal L}_{\zeta_l}
\circ {\Cal L}_{w_l} \circ \cdots \circ
\underline{\Cal L}_{\zeta_1} \circ {\Cal L}_{w_1} (0)),
$$
which are all clearly already analytic, because $h \circ {\Cal
L}_{w_{l+1}} \circ \cdots \circ {\Cal L}_{w_1}(0)$ is analytic by the
induction assumption $({\Cal I}_{2l+1})$ and because
$(\nabla^{\kappa_0} \bar{h}) \circ \underline{\Cal L}_{\zeta_l}
\circ {\Cal L}_{w_l} \circ \cdots \circ
\underline{\Cal L}_{\zeta_1} \circ {\Cal L}_{w_1} (0)$
is analytic by the induction assumption $(*_{2l})$.

Secondly, they also involve terms of the form
$$
\left[ \frac{\partial^{|\alpha|} }{\partial w_{l+1}^{\alpha}}
\frac{\partial^{|\beta|} }{\partial \xi_1^{\beta}} \ \ \left(
\underline{\Upsilon}_{\xi_1} 
\circ {\Cal L}_{w_{l+1}} \circ \cdots 
\circ {\Cal L}_{w_1} (0) \right)\right]_{\xi_1=0}, \ \ \ \ \
|\alpha| \leq |\gamma|+1, \ |\beta| \leq |\delta|,
\tag 8.21
$$
which are obviously analytic.

Thirdly, they involve terms of the form
$$
\left[ \frac{\partial^{|\alpha|} }{\partial w_{l+1}^{\alpha}}
\frac{\partial^{|\beta|} }{\partial \xi_1^{\beta}} \ \left(
(\nabla^{\kappa_0}\bar{h})\circ
\underline{\Upsilon}_{\xi_1} \circ 
\underline{\Cal L}_{\zeta_l} \circ \cdots \circ
{\Cal L}_{w_1}(0)\right) \right]_{\xi_1=0}=
\tag 8.22
$$
$$
=(\underline{\Upsilon}^{\beta} \nabla^{\kappa_0} \bar{h} ) \circ
\underline{\Cal L}_{\zeta_l} \circ \cdots \circ
{\Cal L}_{w_1}(0)= Q_{\alpha, \beta} (\underline{\Cal L}_{\zeta_l}
\circ \cdots \circ {\Cal L}_{w_1} (0), (\nabla^{\kappa_0+|\beta|}
\bar{h}) \circ
\underline{\Cal L}_{\zeta_l} \circ \cdots \circ {\Cal L}_{w_1}(0)),
$$
with $Q_{\alpha, \beta}$ holomorphic in its variables and $Q_{\alpha,
\beta} \equiv 0$ if $|\alpha|>0$, and these terms are firmly analytic
thanks to the induction assumption $(*_{2l})$.

Now, we consider equations (8.18) (or (8.19)) as a $n'\times n'$
matrix equation with the unknowns $({\Cal L}^{\gamma+\1_r^m}
\underline{\Upsilon}^{\delta} h_l) \circ {\Cal L}_{w_{l+1}} \circ
\cdots \circ {\Cal L}_{w_1}(0):=H_k$, $1 \leq k\leq n'$, written in
the form $XH+A=0$. Here, $X$ and $A$ have analytic coefficients, $H$
is formal. But ${\text{\rm det}} \ X$ as a convergent pover series
does not vanish identically, because
$$
{\text{\rm det}} \ X= {\text{\rm det}}\left(
\frac{\partial {\Cal X}_j'}{\partial t_k'}\left(
{\Cal L}_{w_{l+1}} \circ \cdots \circ {\Cal L}_{w_1} (0), \ h \circ
{\Cal L}_{w_{l+1}} \circ \cdots \circ {\Cal L}_{w_1}(0),
\right.\right.
\tag 8.23
$$
$$ \left.\left. (\nabla^{\kappa_0} \bar{h}) \circ
\underline{\Cal L}_{\zeta_l} \circ {\Cal L}_{w_l} \circ \cdots \circ
{\Cal L}_{w_1}(0)\right)_{1\leq j,k\leq n'}\right),
$$
and when we specify the variables $(w_1, \zeta_1, \ldots, w_l,
\zeta_l, w_{l+1})$ in (8.23) above, as $(0, 0, \ldots, 0, 0,
w_{l+1})$, we obtain
$$
{\text{\rm det}} \left( \frac{\partial {\Cal X}_j'}{\partial t_k'}
\left( {\Cal L}_{w_{l+1}} (0), h\circ {\Cal L}_{w_l+1}(0), 
\nabla^{\kappa_0} \bar{h} (0))\right)_{1\leq j,k\leq n'}
\right) \not\equiv_{w_{l+1}} 0 
\ \ \ \hbox{in} \ \C\dl w_{l+1} \dr,
\tag 8.24
$$
by our main assumption (7.4) in Theorem (7.2). Let us denote by $X^T$
the classical adjoint matrix of $X$ that satisfies by definition $X^T
X= X X^T = ({\text{\rm det}} \ X){\text{\rm Id}}_{n'\times n'}$. Then
$X^T$ has convergent power series entries. Also, $({\text{\rm det}} \
X) H= - X^T A$. We deduce that for each $k=1, \ldots, n'$, there
exists converging power series $b$ and $a_k$, $1\leq k\leq n'$,
$a_k\in
\C\{w_1,\ldots, w_{l+1}\}$ with $b={\text{\rm det}} \ X$ such that
$$
b(w_1,\ldots, w_{l+1}) \ {\Cal L}^{\gamma+\1_r^m}
\underline{\Upsilon}^{\delta} h_k({\Cal L}_{w_{l+1}} \circ \cdots
\circ {\Cal L}_{w_1} (0))= a_k(w_1, \ldots, w_{l+1}), \ \ 1\leq k\leq n',
\tag 8.25
$$
with $b\not\equiv 0$. It is then a consequence of the Weierstrass
division theorem that (8.25) implies that there exist $c_j\in
\C\{w_1,\ldots, w_{l+1}\}$ such that $a_k= b\, c_k$. In conclusion,
${\Cal L}^{\gamma+\1_r^m} \underline{\Upsilon}^{\delta} h_k ({\Cal
L}_{w_{l+1}} \circ \cdots \circ {\Cal L}_{w_1}(0)) \in \C\{w_1,\ldots,
w_{l+1}\}$, $1\leq k\leq n'$.

We thus have proved by induction on $\lambda=\vert \gamma \vert +
\vert \delta \vert$ that, for all multiindices $\gamma$, $\delta$,
then $({\Cal L}^{\gamma} \underline{\Upsilon}^{\delta}h ) \circ
\Gamma_{\underline{\Cal L}{\Cal L}}^{2l+1}(w_{l+1}, \zeta_l, w_l, \ldots, 
\zeta_1, w_1) \in \C\{w_1, \zeta_1, \ldots, w_l, \zeta_l, w_{l+1}\}$, 
which proves that $(\nabla^{\kappa} h) \circ \Gamma_{\underline{\Cal
L} {\Cal L}}^{2l+1} (w_{l+1}, \zeta_l, w_l, \ldots,
\zeta_1, w_1)$ 
is convergent for all $\kappa\in \N$.

Finally, it is clear that $(\nabla^{\kappa} \bar{h}) \circ {\Cal
L}_{w_{l+1}} \circ \cdots \circ{\Cal L}_{w_1} (0) = (\nabla^{\kappa}
\bar{h}) \circ \underline{\Cal L}_{\zeta_l} \circ \cdots \circ {\Cal
L}_{w_1}(0) \in \C\{w_1, \ldots, \zeta_l\}$ converges for all
$\kappa$, by the induction assumption $(*_{2l})$.

In conclusion, $\nabla^{\kappa} h^c=(\nabla^{\kappa} h,
\nabla^{\kappa} \bar{h})$ converges on the $(2l+1)$-th Segre chain,
which completes the proof of the implication $[{\Cal I}_{2l+1}$ and
$(*_{2l})]
\Rightarrow (*_{2l+1})$.
\qed
\enddemo

\head \S9. Proof of Theorem 7.2: step II \endhead

\demo{End of proof of Theorem 7.2}
Now, it remains to prove that $(*_{2l+1}) \Rightarrow ({\Cal
I}_{2l+2})$, using $(*_{2l+1})$ that we just have proved.

For that purpose, we start with the conjugate fundamental identities
(8.11) in which we substitute for the variables $(t,\tau)\in {\Cal M}$
the series $\underline{\Cal L}_{\zeta_{l+1}} \circ {\Cal L}_{w_{l+1}}
\circ \cdots \circ
\underline{\Cal L}_{\zeta_1} \circ {\Cal L}_{w_1}(0)$, obtaining
$$
\bar{\Cal X}_j'(\bar{\sigma}\circ 
\underline{\Cal L}_{\zeta_{l+1}} \circ {\Cal L}_{w_{l+1}} \circ \cdots \circ
\underline{\Cal L}_{\zeta_1} \circ {\Cal L}_{w_1} (0), \
\bar{h} \circ 
\underline{\Cal L}_{\zeta_{l+1}} \circ {\Cal L}_{w_{l+1}} \circ \cdots \circ
\underline{\Cal L}_{\zeta_1} \circ {\Cal L}_{w_1} (0), 
\tag 9.1
$$
$$
\ \ \ \ \ \ \ \ \ \ \ 
(\nabla^{\kappa_0} h) \circ {\Cal L}_{w_{l+1}} \circ \cdots \circ
\underline{\Cal L}_{\zeta_1} 
\circ {\Cal L}_{w_1} (0))\equiv 0, \ \ \ \ \
1\leq j\leq n',
$$
since again $(\nabla^{\kappa_0} h) \circ \underline{\Cal
L}_{\zeta_{l+1}}(q(x))\equiv
\nabla^{\kappa_0} h(q(x))$ by Lemma 4.5. 
In eq. (9.1), by $(*_{2l+1})$, all the arguments of $\bar{\Cal X}_j'$
are convergent, except the formal unknowns $\bar{h}_k \circ
\underline{\Cal L}_{\zeta_{l+1}}
\circ \cdots \circ{\Cal L}_{w_1} (0)$. In order to apply
Theorem 1.3.2 to deduce that these unknowns are convergent, it
suffices to observe that the determinant
$$
{\text{\rm det}} \left(
\frac{\partial \bar{\Cal X}_j}{\partial \tau_k''}
\left(\bar{\sigma}\circ
\underline{\Cal L}_{\zeta_{l+1}} 
\circ {\Cal L}_{w_{l+1}} \circ \cdots \circ
{\Cal L}_{w_1} (0), \bar{h} \circ
\underline{\Cal L}_{\zeta_{l+1}} 
\circ {\Cal L}_{w_{l+1}} \circ \cdots \circ
{\Cal L}_{w_1} (0),\right.\right.
\tag 9.2
$$
$$\left.\left. (\nabla^{\kappa_0} h)
\circ {\Cal L}_{w_{l+1}} \circ 
\cdots \circ
\underline{\Cal L}_{\zeta_1}\circ 
{\Cal L}_{w_1} (0)\right)_{1\leq j,k\leq n'}
\right)
$$
does not vanish identically in $\C\dl w_1,\zeta_1, \ldots, w_{l+1},
\zeta_{l+1} \dr$, because, when we specify $(w_1, \zeta_1, \ldots,
w_{l+1}, \zeta_{l+1})$ as $(0, 0, \ldots, 0, \zeta_{l+1})$ in (9.2),
we see that
$$
{\text{\rm det}} \left(
\frac{\partial \bar{\Cal X}_j}{\partial \tau_k'}
\left(\bar{\sigma}\circ
\underline{\Cal L}_{\zeta_{l+1}}(0), \bar{h} \circ 
\underline{\Cal L}_{\zeta_{l+1}}(0),
\nabla^{\kappa_0} h (0)\right)_{1\leq j,l\leq n'}
\right)
\not\equiv 0 \ \ \ \hbox{in} \ \C \dl \zeta_{l+1} \dr,
\tag 9.3
$$
according to our main assumption (3.4).

Finally, it is clear that $h_j \circ \underline{\Cal L}_{\zeta_{l+1}}
\circ \cdots \circ {\Cal L}_{w_1}(0)\equiv h_j \circ {\Cal
L}_{w_{l+1}} \circ \cdots \circ {\Cal L}_{w_1}(0)\in \C\{w_1,\ldots,
w_{l+1}\}$ by the induction assumption $({\Cal I}_{2l+1})$. In
conclusion, $h^c= (h, \bar{h})$ converges on the $(2l+2)$-th Segre
chain, which completes the proof of the implication $(*_{2l+1})
\Rightarrow ({\Cal I}_{2l+2})$.
\qed

The proof of Theorem 7.2 is complete.
\qed\enddemo

\head \S10. Proof of Theorem 7.5: steps I and II \endhead

\demo{Proof of Theorem 7.5}
For the propagation processus, we shall need essentially two lemmas.

\proclaim{Lemma 10.1}
Let $w\in \C^{\mu}$, $z\in \C^d$ and suppose that a power series
$h(w,z)\in \C\dl w, z \dr$ formally satisfies a polynomial equation of
the form
$$
h(w,z)^N + \sum_{1\leq k\leq N} a_k(w,z) \, h(w,z)^{N-k},
\tag 10.2
$$
where $N\in N_*$ and where each power series $a_k(w,z) =
\sum_{\alpha\in
\N^d} z^{\alpha} \, a_{k, \alpha} (w) \in \C \dl w, z \dr$ has all its
derivatives with respect to $z$ at $0$ being convergent power series,
{\it i.e.}
$$
(1/\alpha!) \ \partial_z^{\alpha} a_k(w,0)= a_{k,\alpha} (w) \in
\C\{w\}, \ \ \ \ \ 1\leq k\leq N, \ \forall \ \alpha \in \N^d.
\tag 10.3
$$
Then $h(w,z)=\sum_{\alpha\in \N^d} z^{\alpha} \, h_{\alpha}(w)$ also
satisfies $h_{\alpha}(w)\in \C\{w\}$, $\forall \ \alpha$, {\it i.e.}
$$
(1/\alpha!) \ \partial_z^{\alpha} h(w, 0)= h_{\alpha}(w) \in \C\{w\} \
\
\ \ \ \forall \ \alpha \in \N^d.
\tag 10.4
$$
\endproclaim

We shall observe in a while that the case $d=1$ implies the general
case easily. Let $d=1$. Putting $z=0$ in (10.2), we deduce that
$h(w,0)
\in\C\{w\}$ by applying the following consequence of Artin's theorem.

\proclaim{Lemma 10.5}
Let $w\in \C^{\mu}$, $h\in \C\dl w \dr$, assume that $P(h(w), w)\equiv
0$, where
$$
P(X, w)= \sum_{j=0}^N a_j(w) X^j, \ \ \ \ \ a_j \in \C\{w\}, \ \ \
a_N(w)
\not\equiv 0.
\tag 10.6
$$
Then $h\in \C\{w\}$.
\endproclaim

\demo{Proof}
If $\partial P / \partial X (h(w), w) \equiv 0$, we can of course
replace $P$ by $\partial P / \partial X$. By induction and since
$\partial^N P / \partial X^N = N! \, a_N(w) \not\equiv 0$, we can
assume that $P(h(w), w) \equiv 0$ and $\partial P / P X (h(w), w)
\not\equiv 0$. Finally, applying Theorem 1.3.2, we get $h(w) \in
\C\{w\}$.
\qed\enddemo

\demo{Proof of Lemma 10.1}
Thus $h(w,0)\in \C\{w\}$. Recall that $d=1$. By induction, let us
assume that $h_0(w), \ldots, h_l(w)\in \C\{w\}$ and prove that
$h_{l+1}(w)\in\C\{w\}$. It is clear that if we replace $h$ by
$\tilde{h}(w,z):= h(w,z)-h_0(w)-\cdots- z^lh_l(w)$ in (10.2), we get
immediately that $\tilde{h}$ satisfies a similar polynomial equation
and we have $\tilde{h}= z^{l+1} \tilde{h}_1$. Thus, we can assume
after coming back to the previous notation $h$ that $h=z^{l+1} h_1$
and thus, we must prove that $h_1(w,0)\in \C\{w\}$.

Let us write uniquely $a_k(w,z)= z^{\lambda_k} c_k(w,z)$, where
$\lambda_k\in \N$, $c_k(w,0)\not\equiv 0$ if $a_k(w,z)\not\equiv 0$ or
$\lambda_k=\infty$ if $a_k(w,z)\equiv 0$. Now, in the identity
$$
z^{(l+1)N}\, h_1(w,z)^N + \sum_{1\leq k\leq N} z^{\lambda_k+
(N-k)(l+1)} \, c_k(w,z) \ h_1(w,z)^{N-k}\equiv 0,
\tag 10.7
$$
if we select the term behind $z^{\varpi}$, where
$$
\varpi:= \inf \ \left((l+1)N , 
\inf_{1\leq k\leq N} \lambda_k+(N-k)(l+1)\right) < \infty
\tag 10.8
$$
we shall immediately get that $h_1(w,0)$ satisfies a nontrivial
polynomial equation as in Lemma 10.5. In conclusion, $h_1(w,0)\in
\C\{w\}$ and we are done.

To deduce the general case from the case $d=1$, it suffices to apply
the case $d=1$ to functions $\tilde{h}_c(w, \zeta):= h(w, c\zeta)$,
$w\in \C^{\mu}$, $\zeta\in \C$ for all possible complex lines $\C\ni
\zeta \mapsto (c_1\zeta, \ldots, c_d \zeta) \in \C^d$, $(c_1,\ldots,
c_d)\in \C^d$.
\qed\enddemo
\enddemo

\demo{End of proof of Theorem 7.5}
First, taking $\tau=0$ and $\rho(t, 0)=0$ in (7.7), we deduce that $h
\circ \Gamma_{\underline{\Cal L}{\Cal L}}^1 (w_1) \in
\C\{w_1\}$ by a straightforward application of Lemma 10.1 above.

As in the proof of Theorem 7.2, it suffices to establish that $[({\Cal
I}_k)$ and $(*_{k-1})] \Rightarrow (*_k)$ and that $(*_k)\Rightarrow
({\Cal I}_{k+1})$ in case $k=2l+1$ is odd.

Firstly, for the first step, we substitute the series
$$
\underline{\Upsilon}_{\xi_1} \circ 
\Gamma_{\underline{\Cal L}{\Cal L}}^{2l+1}
(w_{l+1}, \zeta_l, w_l, \ldots, \zeta_1, w_1)\in
\C\{w_1, \zeta_1, \ldots, w_l, \zeta_l, w_{l+1}, \xi_1\}^{2n}
\tag 10.9
$$
for the variables $(t,\tau)\in {\Cal M}$ in the fundamental identity
(7.7) and get
$$
P_j'(\underline{\Upsilon}_{\xi_1} \circ {\Cal L}_{w_{l+1}} \circ
\cdots \circ
\underline{\Cal L}_{\zeta_1}\circ {\Cal L}_{w_1} (0), \
h_j\circ \underline{\Upsilon}_{\xi_1}
\circ {\Cal L}_{w_{l+1}} \circ \underline{\Cal L}_{\zeta_l} \circ 
{\Cal L}_{w_l}\circ \cdots \circ
\underline{\Cal L}_{\zeta_1}\circ {\Cal L}_{w_1} (0),
\tag 10.10
$$
$$ \ \ \ \ \ \ \ \ \ \ \ \ \ \ \ (\nabla^{\kappa_0} \bar{h}) \circ
\underline{\Upsilon}_{\xi_1} \circ 
\underline{\Cal L}_{\zeta_l} \circ 
 {\Cal L}_{w_{l}} \circ \cdots \circ
\underline{\Cal L}_{\zeta_1} \circ {\Cal L}_{w_1} (0)\equiv 0, \ \ \ \ 
1\leq j\leq n',
$$
after the simplification $(\nabla^{\kappa_0} \bar{h})
\circ \underline{\Upsilon}_{\xi_1} \circ {\Cal L}_{w_{l+1}} (q(x))\equiv
(\nabla^{\kappa_0} \bar{h})
\circ \underline{\Upsilon}_{\xi_1} (q(x))$.
Since all derivatives
$$
\left[\frac{\partial^{|\beta|} }{\partial w_{l+1}^{\beta}}
\frac{\partial^{|\alpha|} }{\partial \xi_1^{\alpha}} \left(
(\nabla^{\kappa_0} \bar{h}) \circ
\underline{\Upsilon}_{\xi_1} \circ 
\underline{\Cal L}_{\zeta_l} \circ 
 {\Cal L}_{w_{l}} \circ \cdots \circ
\underline{\Cal L}_{\zeta_1} \circ {\Cal L}_{w_1} (0)\right)\right]_{\xi_1=0}
\tag 10.11
$$
are converging, by assumption $(*_{2l})$, we can apply Lemma 10.1 to
deduce that the derivatives
$$
\left[
\frac{\partial^{|\beta|} }{\partial w_{l+1}^{\beta}}
\frac{\partial^{|\alpha|} }{\partial \xi_1^{\alpha}} \left(
h_j \circ
\underline{\Upsilon}_{\xi_1} \circ 
{\Cal L}_{w_{l+1}} \circ
\underline{\Cal L}_{\zeta_l} \circ 
 {\Cal L}_{w_{l}} \circ \cdots \circ
\underline{\Cal L}_{\zeta_1} \circ {\Cal L}_{w_1} (0)\right)
\right]_{\xi_1=0}
\tag 10.12
$$ 
are in $\C\{w_1,\ldots, \zeta_l, w_{l+1}\}$, for all $\alpha\in \N^d$,
$\beta\in \N^m$: this is $(*_{2l+1})$.

Secondly, for the second step, considering the conjugate identities
(8.11), we get
$$
\bar{P_j}(\bar\sigma\circ
\underline{\Cal L}_{\zeta_{l+1}} \circ {\Cal L}_{w_{l+1}} \circ \cdots \circ
\underline{\Cal L}_{\zeta_1} \circ {\Cal L}_{w_1} (0), \
\bar{h}_j \circ \underline{\Cal L}_{\zeta_{l+1}} 
\circ {\Cal L}_{w_{l+1}} \circ \cdots \circ
\underline{\Cal L}_{\zeta_1} \circ {\Cal L}_{w_1} (0), 
\tag 10.13
$$
$$
\ \ \ \ \ \ \ \ \ \ \ \ \
(\nabla^{\kappa_0} h) \circ {\Cal L}_{w_{l+1}}
\circ \cdots \circ \underline{\Cal L}_{\zeta_1}
\circ {\Cal L}_{w_1} (0)) \equiv 0, \ \ \ \ 
1\leq j\leq n'.
$$
In (10.13) above, the formal series in $(\nabla^{\kappa_0}h) \circ
{\Cal L}_{w_{l+1}} \circ \cdots \circ
\underline{\Cal L}_{\zeta_1} \circ {\Cal L}_{w_1} (0)$ 
are converging, by assumtion $(*_{2l+1})$. Thus, applying Lemma 10.5,
we deduce that $\bar{h} \circ \Gamma_{\underline{\Cal L}{\Cal
L}}^{2l+1} (w_1,
\zeta_1, \ldots, w_{l+1}, \zeta_{l+1})\in \C\{w_1, \zeta_1, 
\ldots, w_{l+1}, \zeta_{l+1}\}^{n'}$. 

Finally, it is clear that $h\circ \underline{\Cal L}_{\zeta_{l+1}}
\circ {\Cal L}_{w_{l+1}} \circ \cdots \circ
\underline{\Cal L}_{\zeta_1}\circ {\Cal L}_{w_1} (0) = 
h\circ \Cal L_{w_{l+1}} \circ \cdots \circ
\underline{\Cal L}_{\zeta_1} \circ {\Cal L}_{w_1} (0) \in
\C\{w_1, \zeta_1, \ldots, 
w_{l+1}\}^{n'}$ by $({\Cal I}_{2l+1})$. In conclusion, $h^c= (h,
\bar{h})$ converges on the $(2l+2)$-th Segre chain.

This completes the proof of Theorem 7.5.
\qed\enddemo

\head \S11 Simplifications in the S-solvable case \endhead

In the case where $h$ is S-solvable at $0$, since we can solve $h$ in
terms of $\bar h$ and its jets $\nabla^\kappa \bar h$ by the usual
analytic implicit function theorem, which takes place instead of
Artin's approximation theorem, then the propagation process can be
highly simplified.

\demo{Proof of Theorem 1.2.1 (i)}
Indeed, according to Lemma 7.12, we can write
$$
h_j'(\zeta, \xi)\equiv \Psi_j(w,z,\zeta,\xi,
\nabla^{\kappa_0}\bar{h}(\tau)), \ \ \ \ \
1\leq j\leq n',
\tag 11.1
$$
for holomorphic $\Psi_j$ in terms of their arguments, when $(t,\tau)$
satisfy $\rho(t,\tau)=0$. Furthermore, it is possible ({\it cf.}
[BER97]) to obtain after applying the $d$-vector fields $\Upsilon$ and
$\underline{\Upsilon}$ to the above identity (11.1) to obtains for all
$\kappa \in \N$ the existence of holomorphic $\Psi^{\kappa}$ and
$\underline{\Psi}^{\kappa}$ such that
$$
\nabla^{\kappa} \bar{h}(\zeta,\xi)\equiv_{\zeta,w,\xi}
\left[\Psi^{\kappa}(w,z, \zeta, \xi, 
\nabla^{\kappa+\kappa_0} {h} (w,z))\right]_{z:=
\xi+i\bar{\Theta}(\zeta,w,\xi)}
\tag 11.2
$$
$$\nabla^{\kappa} h(w,z) \equiv_{w,\zeta,z}
\left[\underline{\Psi}^{\kappa} 
(\zeta, \xi, w, z, \nabla^{\kappa+\kappa_0}
\bar{h}(\zeta, \xi))\right]_{\xi:=z-i\Theta(w,\zeta,z)}.
$$
But starting with the more intrinsic (equivalent) writing of (11.2) as
follows
$$ (\nabla^{\kappa}\bar{h}) \circ
\underline{\Cal L}_{\zeta} (p)
=\underline{\Psi}^{\kappa} (\underline{\Cal L}_{\zeta}(p),
\nabla^{\kappa+\kappa_0} h(p)),
\tag 11.3
$$
$$\ \ \ \ \ \ \ \ \ (\nabla^{\kappa} h) \circ {\Cal
L}_w(p)=\Psi^{\kappa} ({\Cal L}_w(p), \nabla^{\kappa+\kappa_0} \bar{h}
(p)),
\ \ \ \ p=p^c\in {\Cal M} \cap 
\underline{\Lambda}, \ \ \ \kappa\in \N,
$$
we can make immediate iterations of identities (11.3), for instance
twice:
$$
(\nabla^{\kappa} h) \circ {\Cal L}_{w_1} \circ
\underline{\Cal L}_{\zeta_1} (p) =\Psi^{\kappa}(
{\Cal L}_{w_1} \circ \underline{\Cal L}_{\zeta_1} (p),
\nabla^{\kappa+\kappa_0} \bar{h}
\circ \underline{\Cal L}_{\zeta_1} (p)=
\tag 11.4
$$
$$ 
\ \ \ \ \ \ \ \ \ \ \ \ 
\Psi^{\kappa}({\Cal L}_{w_1} \circ 
\underline{\Cal L}_{\zeta_1} (p), 
\Psi^{\kappa+\kappa_0}(\underline{\Cal L}_{\zeta_1}(p),
\nabla^{\kappa+2\kappa_0} h(p))),
$$
$$ \ \ \ \ \ \ \ \ \ \ \ \ =
\psi_2^{\kappa} (w_1,\zeta_1, \nabla^{\kappa+2\kappa_0}
h(p))\in \C\{w_1,\zeta_1\}
$$
so that to perform the propagation process, we only have to iterate
({\it i.e.} replace them into themselves) the identities (11.2), and
thus, we immediately obtain
$$ 
(\nabla^{\kappa} \bar{h})\circ
\underline{\Cal L}_{\zeta_{\mu_p}} 
\circ {\Cal L}_{w_{\mu_p}} \circ \cdots \circ
\underline{\Cal L}_{\zeta_1} \circ {\Cal L}_{w_1} (p) \equiv
\underline{\psi}_{2\mu_p}^{\kappa} (\zeta_1,w_1,...,
\zeta_{\mu_p}, w_{\mu_p},
\nabla^{\kappa+2\mu_p\kappa_0} \bar{h} (p)),
\tag 11.5
$$
$$ \ \ \ \ \ (\nabla^{\kappa} h )\circ {\Cal L}_{w_{\mu_p}}
\circ \underline{\Cal L}_{\zeta_{\mu_p}} 
\circ \cdots \circ{\Cal L}_{w_1} \circ 
\underline{\Cal L}_{\zeta_1} (p)\equiv
\psi_{2\mu_p}^{\kappa}(
w_1,\zeta_1,...,w_{\mu_p}, \zeta_{\mu_p},
\nabla^{\kappa+2\mu_p \kappa_0} h(p)),
$$
which, together with Theorem 3.2.2, completes the proof of Theorem
1.2.1 (i).
\qed
\enddemo

\heading \S12. Analyticity of formal solutions \endheading

In this very short paragraph, we prove Theorem 1.3.2, which has been
applied several times during the proof of Theorem 1.2.1. By ${\Cal
V}_X(p)$, we mean a small open neighborhood, equivalent to a polydisc,
of the point $p$ in the complex manifold $X$ equivalent to a polydisc.

\demo{Proof of Theorem 1.3.2} 
Fix $j_1,\ldots,j_m$ with
$$
\text{\rm det} \left(
\frac{\partial R_{j_k}}{\partial y_l}(w,\hat{g}(w))
\right)_{1\leq k,l \leq m}\not\equiv_w 0 \ \ \text{\rm in} \ \ \C\dl w\dr
\tag 12.1
$$
and simply write $R_1,\ldots,R_m$. We show that the solution
$(\hat{g}_1(w),\ldots,\hat{g}_m(w))$ of the $m$ equations
$R_1(w,\hat{g}(w))\equiv_w 0,\ldots, R_m(w,\hat{g}(w))\equiv_w 0$ is
itself converging.

Let $\Gamma$ be the germ at $0$ of the complex analytic set
$$ 
\Gamma=\{(w,y)\in {\Cal V}_{\C^n}(0)\times
{\Cal V}_{\C^m}(0) \: R_1(w,y)=0,\ldots,R_m(w,y)=0\}.
\tag 12.2
$$
For each $N\in \N$, there exists an analytic solution $g_N(w)$, {\it
i.e.} $R(w,g_N(w))\equiv_w 0$, with $g-g_N = 0 \, (\text{\rm mod} \
\frak{m}(w)^N)$. Consider the graph of $g_N$:
$$ 
\text{\rm gr} (g_N) =\{(w,y) \in {\Cal V}_{\C^n}(0) 
\times {\Cal V}_{\C^n}(0) \: 
y=g_N(w)\}.
\tag 12.3
$$
For $N$ large enough, say $N\geq N_0$, $N_0\in \N_*$, by (12.1),
$$
\text{\rm det} \left(
\frac{\partial R_k}{\partial y_l}(w,g_N(w))
\right)_{1\leq k,l \leq m}\not\equiv_w 0,
\tag 12.4
$$
so that there exists $w_0$ close to $0$ in $\C^n$ such that
$$ 
\text{\rm det} \left(
\frac{\partial R_k}{\partial y_l}(w_0,g_N(w_0))
\right)_{1\leq k,l \leq m} \neq 0.
\tag 12.5
$$
By (12.5), a neighborhood of $(w_0, g_N(w_0))$ in $\Gamma$ is of
dimension $n$. Therefore,
$$
\text{\rm gr}(g_N) \cap 
({\Cal V}_{\C^n}(w_0) \times {\Cal V}_{\C^m} (g_N(w_0)))\equiv
\Gamma \cap 
({\Cal V}_{\C^n}(w_0) \times {\Cal V}_{\C^m} (g_N(w_0))).
\tag 12.6
$$
As a consequence, {\it for each $N\geq N_0$, ${\text{\rm gr}}(g_N)$ is
an irreducible component of $\Gamma$.} Since $\Gamma$ has finitely
many irreducible components, a subsequence of $g_N$, {\it i.e.} of
${\text{\rm gr}}(g_N)$, is constant. Since $g_N \to \hat{g}$ in the
Krull topology, $\hat{g}$ is convergent.
\qed\enddemo

\remark{Remark}
In fact, more is {\it a posteriori} true above: it is clear then that
there exists an integer $N_0$, depending only on the $R_k$'s, $1\leq
k\leq m$, such that for all $N, N'\geq N_0$, $g_N=g_{N'}$ near $0$
(the component stabilizes). More generally, we observe:
\endremark

\proclaim{Lemma 12.7}
Let $R_1(w,y),\ldots,R_m(w,y)\in\C\{w,y\}$ and assume that there
exists $g^1(w)\in \C\dl w\dr^m$ satisfying
$$
R_j(w,g^1(w))\equiv_w 0 \ \ \ \hbox{and}
\ \ \
\text{\rm det} \left(
\frac{\partial R_{k}}{\partial y_l}(w,g^1(w))
\right)_{1\leq k,l \leq m}\not\equiv_w 0 \ \ \text{\rm in} \ \ \C\dl w\dr
\tag 12.8
$$
Then there exists a positive integer
$\underline{\nu}=\underline{\nu}(R)$ such that if $g^2\in \C\dl w\dr$
satisfies
$$
\aligned
& R_k(g^2(w), w)\equiv 0 \ \ \ \ \ \hbox{in} \ \ \C\dl w\dr, \ \ \ \ \
\forall \ 1\leq k\leq m\\ &
\partial_w^\alpha g^1(0) = \partial_w^\alpha g^2(0) \ \ \ \ \ \forall \
\v \alpha \v \leq \underline{\nu}(R),
\endaligned
\tag 12.9
$$
then $g^2(w)\equiv g^1(w)$.
\endproclaim

\demo{Proof}
By the above proof of Theorem 1.3.2, $\hbox{gr}(g^1)$ is then a fixed
irreducible component $\Gamma^1$ of the complex analytic set $\Gamma$
defined by (12.2). It is also clear that there exists a sufficiently
large integer $\underline{\nu}$ depending only on $R$ such that
$\partial_w^\alpha g^1(0) = \partial_w^\alpha g^2(0), \, \forall \
\v \alpha \v \leq \underline{\nu}$ implies that 
$\text{\rm det} \left(
\frac{\partial R_{k}}{\partial y_l}(w,g^2(w))
\right)_{1\leq k,l \leq m}\not\equiv_w 0$. In this case, $\hbox{gr}(g^2)$
occurs to be also an irreducible component $\Gamma^2$ of $\Gamma$.
Furthermore, $\Gamma^1=\Gamma^2$ if $\underline{\nu}=
\underline{\nu}(R)$ is large enough. Finally, $g^1(w)\equiv g^2(w)$.
\qed
\enddemo

In particular, a direct corollary is as follows ({\it cf.} [BER99],
Lemma 4.3):

\proclaim{Corollary 12.10}
Let $P(X,w)$ be of the form
$$
P(X,w)=\sum_{j=0}^N a_j(w)X^j, \ \ \ \ \ a_j\in \C\{w\}, \ \
a_N(w)\not\equiv 0.
\tag 12.11
$$
Then there exists a positive integer $\underline{\nu}(P)$ such that if
$h^1, h^2\in
\C\dl w\dr$ satisfy
$$
\aligned
& P(h^1(w), w) \equiv P(h^2(w),w)\equiv 0 \ \ \ \ \ \hbox{in} \ \C\dl
w\dr\\ &
\partial_w^\alpha h^1(0) = \partial_w^\alpha h^2(0) \ \ \ \ \ \forall \
\v \alpha \v \leq \underline{\nu}(P),
\endaligned
\tag 12.13
$$
then $h^1(w)\equiv h^2(w)$.
\endproclaim

\demo{Proof}
The only thing we have to check is that we can assume in addition that
$\partial_XP(h^1(w),w)\not\equiv 0$. We leave this to the reader.
\qed\enddemo

\heading \S13. Examples\endheading

\example{Example 13.1}
{\it An S-solvable or S-finite formal map $h\: (M,p)\to_{\Cal F} (M',
p')$ such that the formal generic rank of $h\: (S_{\bar{p}}, p)
\to_{\Cal F} (S_{\bar{p}'}', p')$ is strictly less than $m'$}. Take
$(z_1, z_3)\mapsto (z_1,0,z_3)$, $\C_{z_1, z_3}^2 \to
\C_{z_1',z_2',z_3'}^3$, $M=\{ z_3=\bar{z}_3+iz_1\bar{z}_1\}$ and
$M'=\{z_3'=\bar{z}_3'+i z_1'\bar{z}_1'+i {z'}_1^a {\bar{z'}}_2^a+
i{\bar{z'}}_1^a{z'}_2^a$, $a\in N_*$. If $a=1$, then $h$ is S-solvable
at the origin. If $a\geq 2$, then $h$ is S-finite but not S-solvable
at the origin (of course, $M'$ is essentially finite at the origin.)
\endexample

\example{Example 13.2}
{\it An S-finite but not S-nondegenerate formal map}. Take $(z_1,
z_3)\mapsto (z_1,0,z_3)$, $\C_{z_1, z_3}^2 \to \C_{z_1',z_2',z_3'}^3$,
$M=\{ z_3=\bar{z}_3+iz_1^2\bar{z}_1^2\}$ and $M'=\{ z_3'=\bar{z}_3' +i
{z'}_1^2 {\bar{z'}}_1^2+ i{z}_1' {\bar{z'}}_2^2 + i{\bar{z'}}_1
{z'}_2^2\}$.
\endexample

\example{Example 13.3}
{\it An S-nondegenerate but not S-finite formal map}. Perharps the
simplest example is the identity map of $M=\{z_3=\bar{z}_3 + i z_1
\bar{z}_1 (1+ z_2\bar{z}_2)\}$ ({\it cf.} [MM2]).
\endexample

\example{Example 13.4}
{\it A holomorphically nondegenerate but not S-nondegenerate real
hypersurface}. $M \: \, y_3=\vert z_1 \vert^2 \vert 1+ z_1
\bar{z}_2\vert^2 (1+\hbox{Re} (z_1\bar{z}_2))^{-1} -x_3\,
\hbox{Im} (z_1\bar{z}_2) 
(1 +\hbox{Re} (z_1\bar{z}_2))^{-1}$. This seems to be the simplest
example of such ({\it cf.} [BER99] [Mer99c]).
\endexample

\example{Example 13.5}
{\it A class of S-finite or S-nondegenerate real hypersurfaces in
$\C^n$.} Here are the two most naive examples: The hypersurface $y=\v
w_1\v^{2r_1}+\cdots+\v w_{n-1}\v^{2r_{n-1}}+ xh(w,\bar w, x)$ where
$h$ is any analytic remainder in normal form and where
$r_1>0,\ldots,r_{n-1}>0$, is essentially finite (S-finite). The
hypersurface $y=\sum_{k=1}^{\mu} \prod_{j=1}^{n-1} \v w_j
\v^{2r_{j,k}} (1+xh(w,\bar w,x))$, $\mu \geq n-1$, where $h$ is any
analytic remainder and where the exponents $r_{j,k}$ are subject to
the condition that the generic rank of the Jacobian matrix
$\left({\partial\prod_{j=1}^{n-1} w_j^{r_{j,k}}\over
\partial w_l}\right)_{k,l}$ equals $n$ and all the 
multiindinces $\beta_k:=(r_{j,k})_{1\leq j\leq n-1}\in \N_*^{n-1}$ are
pairwise distinct.

\endexample
\head \S14. A uniqueness principle \endhead

The purpose of this paragraph is to establish the following uniqueness
principle, first derived in the S-solvable case in [BER97] ({\it cf.}
also [BER99], Theorem 2.5):

\proclaim{Theorem 14.1}
Let $h^1\: (M,p)\to_{\Cal F} (M',p')$ be a real analytic holomorphic
mapping between real analytic CR generic manifolds, assume that $M$ is
minimal at $p$, and assume that $h^1$ is either S-solvable, or
S-finite, or S-nondegenerate. Then there exists an integer
$\underline{\kappa}\in \N_*$ with the following property. If $h^2\:
(M,p)\to_{\Cal F} (M',p')$ is a formal holomorphic mapping sending
$(M,p)$ into $(M',p')$, and if
$$
\partial_t^\alpha h^2(p)=\partial_t^\alpha h^1(p) \ \ \ \ \ 
\forall \ \v \alpha \v \leq \underline{\kappa},
\tag 14.2
$$
then $h^1(t-p)\equiv h^2(t-p)$ in $\C\dl t-p\dr$.
\endproclaim

First, here an immediate consequence in the invertible case:

\proclaim{Corollary 14.3}
Let $M$ be a real analytic generic submanifold through the origin in
$\C^n$. Assume that $(M,0)$ is minimal and S-solvable, or essentially
finite, or S-nondegenerate at 0. Then there exists an integer
$\underline{\kappa}$ with the following property. If $M'$ is a real
analytic generic submanifold through the origin in $\C^n$ of the same
dimension as $M$ and if $h^1, h^2\: (\C^n, 0)\to_{\Cal F} (\C^n,0)$
are formal invertible mappings sending $(M,0)$ into $(M',0)$ which
satisfy
$$
\partial_t^\alpha h^1(t) = \partial_t^\alpha h^2(0) \ \ \ \ \
\forall \ \v \alpha \v \leq \underline{\kappa},
\tag 14.4
$$
then $h^1(t)\equiv h^2(t)$ in $\C\dl t\dr$.
\endproclaim

\demo{Proof}
We claim that it suffices to take $\underline{\kappa}$ to be the
integer given by Theorem 14.1 with $M=M'$ and $h^1=\hbox{Id}$, {\it
i.e.} $h^1(t)\equiv t$. To see this, let $M'$, $h^1$, $h^2$ be as in
Corollary 14.3 and observe that if (14.4) holds, then
$\partial_t^\alpha ((h^1)^{-1} \circ h^2)(0)=0$ for all $\v \alpha \v
\leq
\underline{\kappa}$. By Theorem 14.1, and the choice of $\underline{\kappa}$, 
we deduce that $((h^1)^{-1}\circ h^2)(t)\equiv t$ and hence, the
conclusion of Corollary 14.3.
\qed
\enddemo

\demo{Proof of Theorem 14.1}
To prove Theorem 14.1, we follow the steps of the proof of Theorem
1.3.2 thoroughly. It suffices to treat S-finite and S-nondegenerate
maps parallely (S-solvable maps being S-finite). We shall concentrate
on the S-nondegenerate case only (to treat the S-finite case, use
Corollary 12.10 instead of Lemma 12.7). In what follows, $\Cal
F(M,M')$ we denote the set of all formal mappings $(\C^n,0)\to_{\Cal
F} (\C^{n'},0)$ that send $(M,0)$ into $(M',0)$. We consider the
following property for $k\in \N$ and $\kappa \in \N$.

\smallskip
\noindent
{\it $(**)_{k,\kappa}$ There exists $K(k,\kappa)\in \N$ such that for
any $h^2\in \Cal F(M,M')$ with
$$
\partial_t^\alpha h^1(0) = \partial_t^\alpha h^2(0), \ \ \ \ \ 
\forall \ \v \alpha \v \leq K(k,\kappa),
\tag 14.5
$$
the following holds}
$$
(\nabla^\kappa h^{1c})\circ \Gamma_{\underline{\Cal L} {\Cal
L}}^k(w_1,\ldots,w_k)\equiv (\nabla^\kappa h^{2c})\circ
\Gamma_{\underline{\Cal L} {\Cal L}}^k(w_1,\ldots,w_k).
\tag 14.6
$$

Observe that $(**)_{0,\kappa}$ holds with $K(0,\kappa)=\kappa$, since
$\Gamma_{\underline{\Cal L} {\Cal L}}^0\equiv 0$ is the constant null
mapping. We shall prove that $(**)_{k,\kappa}$ holds for all $k$ and
$\kappa$ by double induction on $k$ and $\kappa$. First, let us assume
that $(**)_{k' ,\kappa'}$ holds for all $0\leq k'\leq 2l+1$ and all
$\kappa'$ and prove that $(**)_{2l+2, 0}$ holds (as in \S8-9, the
even-to-odd induction is similar). Precisely, we must show the
existence of the integer $K(2l+2,0)$ in $(**)_{2l+2,0}$. Let
$\bar{\Cal X}_j'(\tau,t,\tau',
\nabla^{\kappa_0})$, $1\leq j\leq n'$, be the (conjugate, {\it cf.} eq.
(8.11)) analytic functions given in the assumptions of Theorem 7.2.
Pick an integer $\tilde{K}$ and consider the formal mappings $h^2\in
\Cal F(M,M')$ satisfying
$$
\partial_t^\alpha h^1(0) = \partial_t^\alpha h^2(0), \ \ \ \ \ 
\forall \ \v \alpha \v \leq \tilde K.
\tag 14.7
$$
By {\it a priori} \, requiring $\tilde K\geq K(2l+1,\kappa_0)$, we may
assume that any $h^2\in \Cal F(M,M')$ that satisfies (14.7) above also
satisfies the identities (8.11) {\it on the $2l+2$-th Segre chain},
with $\bar h^1$ replaced by $\bar h^2$, and the same jets
$\nabla^{\kappa_0} h^1$, {\it i.e.} that
$$
\bar{\Cal X}_j'(\bar{\sigma}\circ 
\underline{\Cal L}_{\zeta_{l+1}} \circ {\Cal L}_{w_{l+1}} \circ \cdots \circ
\underline{\Cal L}_{\zeta_1} \circ {\Cal L}_{w_1} (0), \
\bar{h}^2 \circ 
\underline{\Cal L}_{\zeta_{l+1}} \circ {\Cal L}_{w_{l+1}} \circ \cdots \circ
\underline{\Cal L}_{\zeta_1} \circ {\Cal L}_{w_1} (0), 
\tag 14.8
$$
$$
\ \ \ \ \ \ \ \ \ \ \ 
(\nabla^{\kappa_0} h^1) \circ {\Cal L}_{w_{l+1}} \circ \cdots \circ
\underline{\Cal L}_{\zeta_1} 
\circ {\Cal L}_{w_1} (0))\equiv 0, \ \ \ \ \
1\leq j\leq n',
$$
Hence, if $\underline{\nu}^{2l+2}$ is the integer given by Lemma 12.7
for equations (14.8), and if we choose $K(2l+2,0)=\max (\tilde K,
\underline{\nu}^{2l+2})$, then the identity $h^1\circ
\Gamma_{\underline{\Cal L} {\Cal L}}^{2l+2}= h^2\circ
\Gamma_{\underline{\Cal L} {\Cal L}}^{2l+2}$ follows from Lemma
12.7. The property $(**)_{2l+2,0}$ is proved.

We now fix $k=2l+1$ and an integer $\kappa$. We complete the induction
by assuming that $(**)_{k',\kappa'}$ for all pairs $(k', \kappa')$
satisfying either $0\leq k' < 2l+1$ or $k' = 2l+1$ and
$\kappa'\leq\kappa$ (the case $k=2l$ is similar). We shall prove
$(**)_{2l+1,\kappa+1}$. Now, consider those $h^2\in \Cal F(M,M')$
satisfying (14.7) with $\tilde K \geq \max \{
\tilde K(2l,\kappa_0+\kappa+1), \tilde K(2l+1,\kappa)\}$. Using the
induction property, we see that, as above, the derivatives $\Cal
L^{\gamma+\1_r^m} \underline{\Upsilon}^\delta h^2$ and $\Cal L^\gamma
\underline{\Upsilon}^{\delta+\1_s^d} h^2$ of such an $h^2$ satisfy the
same identities (8.18) and (8.19) as $\Cal L^{\gamma+\1_r^m}
\underline{\Upsilon}^\delta h^1$ and $\Cal L^\gamma
\underline{\Upsilon}^{\delta+\1_s^d} h^1$ on the $2l+1$-th Segre chain
(with the same $A_{j,\gamma,\delta,\1_r^m}$ and
$B_{j,\gamma,\delta,\1_s^d}$ for both), for all $\gamma, \delta$ with
$\v \gamma \v + \v \delta \v \leq \kappa$. It is then clear that after
taking the inverse of the matrix (8.23) as in the end of \S8, we
immediately get the agreement of $\Cal L^{\gamma+\1_r^m}
\underline{\Upsilon}^\delta h^2$ and $\Cal L^\gamma
\underline{\Upsilon}^{\delta+\1_s^d} h^2$ with $\Cal
L^{\gamma+\1_r^m} \underline{\Upsilon}^\delta h^1$ and $\Cal L^\gamma
\underline{\Upsilon}^{\delta+\1_s^d} h^1$ on the $2l+1$-th Segre
chain. This completes the induction and proves $(**)_{k,\kappa}$ for
all $k$ and $\kappa$.

To complete the proof of Theorem 14.1, it suffices to remember the
full rank property of $\Gamma_{\underline{\Cal L} {\Cal L}}^{2\mu_0}$
and to apply it to $(**)_{2\mu_0,0}$: $h^1\circ
\Gamma_{\underline{\Cal L} {\Cal L}}^{2\mu_0}=h^2\circ
\Gamma_{\underline{\Cal L} {\Cal L}}^{2\mu_0}$, whence $h^1=h^2$.
This completes the proof of Theorem 14.1.
\qed
\enddemo

\head \S15. Open problems \endhead

First, we wonder whether an analog of Artin's theorem for CR maps
holds true:

\proclaim{Problem 15.1}
Let $h\: (X,p)\to_{\Cal F} (X',p')$ be a formal holomorphic map
between two germs of real analytic sets in $(\C^n,p)$ and
$(\C^{n'},p')$ respectively. Prove $($or disprove$)$ that for every
$N\in
\N$ there exists a convergent mapping $h_N\: (X,p)\to (X',p')$
such that $h(t)\equiv h_N(t) \ (\text{\rm mod} \, \frak{m} (t)^N)$.
\endproclaim

Our analysis leaves open at least four more specialized questions:

\proclaim{Problem 15.2}
Prove that a formal holomorphic mapping $h\: (M,p)\to_{\Cal F}
(M',p')$ between real analytic CR manifolds, such that $(M,p)$ is
minimal and $(M',p')$ does not contain complex analytic sets of
positive dimension through $p'$, is convergent.
\endproclaim

\remark{Remark}
The hypersurface case is treated in [BER99]. It seems that a plausible
adaptation of double reflection of jets as in [Z98] [D] [Mer99b] would
yield the result modulo some unavoidable technicalities.
\endremark

A CR manifold $M$ is called {\it transversally nondegenerate} if it is
not biholomorphically equivalent, locally around a generic point, to a
product $\underline{M} \times I$ of a CR manifold $\underline{M}
\subset \C^{n-1}$ with a real segment $I\subset \R$ ([Mer99a]). A
hypersurface $M$ is transversally nondegenerate if and only if it is
Levi-nonflat.

\proclaim{Problem 15.3}
Prove $($or disprove$)$ that any formal invertible CR map between
holomorphically nondegenerate transversally nondegenerate real
analytic CR manifolds is convergent. Study formal invertible self maps
of nonminimal hypersurfaces.
\endproclaim

\remark{Remark}
Even in codimension one, this problem is new and widely open. In
codimension two, perharps the simplest example to study would be self
maps of $(M,p)=(M',p')=(N,0)$ where $N: z_3=\bar z_3, \, z_2=\bar
z_2+iz_1\bar z_1+i(z_1\bar z_1)^2\bar z_3$.
\endremark

\proclaim{Problem 15.4}
Prove that the reflection mapping associated with a formal CR map from
$(M,p)$ ${\Cal C}^\omega$ minimal into any ${\Cal C}^\omega$ $(M',p')$
is always convergent $($without any nondegeneracy condition on
$(M',p'))$.
\endproclaim

\proclaim{Problem 15.5}
Find purely algebraic proofs of the convergence results in the
circumstance where both $(M,p)$ and $(M',p')$ are algebraic.
\endproclaim

\enddocument
\end